\documentclass{article}
\usepackage{amsmath,amssymb,amsthm}
\usepackage{thmtools}
\usepackage{hyperref}
\usepackage{mathtools}
\usepackage{xcolor}
\usepackage[nameinlink]{cleveref}
\usepackage{parskip}
\usepackage{enumitem}
\usepackage{comment}
\usepackage[
backend=biber,
style=alphabetic,
maxbibnames=99,
maxcitenames=99
]{biblatex}
\bibliography{references}
\usepackage[scale=0.75]{geometry}

\usepackage{titlesec}
\titleformat{\paragraph}[block]{\bfseries\filcenter}{}{0pt}{}

\newcommand{\R}{\mathbb R}
\newcommand{\Z}{\mathbb Z}
\newcommand{\Q}{\mathbb Q}
\newcommand{\N}{\mathbb N}

\newcommand{\Legendre}[2]{\genfrac{(}{)}{}{}{#1}{#2}}

\newtheoremstyle{plain}
    {\medskipamount}
    {}
    {\itshape} 
    {}    
    {\bfseries}
    {.}     
    { }   
    {}       
\theoremstyle{plain}

\newtheorem{theorem}{Theorem}[section]
\newtheorem{proposition}[theorem]{Proposition}
\newtheorem{definition}[theorem]{Definition}
\newtheorem{remark}[theorem]{Remark}

\newtheorem{corollary}[theorem]{Corollary}
\newtheorem{lemma}[theorem]{Lemma}
\newtheorem{conjecture}[theorem]{Conjecture}

\DeclareMathOperator{\sqf}{sqf}
\DeclareMathOperator{\sqp}{sqp}

\newcommand{\stirlingii}[2]{\genfrac{\{}{\}}{0pt}{}{#1}{#2}}
\renewcommand{\Re}{\operatorname{Re}}
\renewcommand{\Im}{\operatorname{Im}}

\title{On the Distribution of Prime Factors of $\lfloor \alpha n^2 \rfloor$}
\author{Cati Alexandru, Johannes Borgwardt, Benjamin Faltin, Giorgos Petritakis}
\date{}

\begin{document}

\maketitle
\begin{abstract}
    Let $\alpha > 0$ be irrational. We show that $\omega(\lfloor\alpha n^2\rfloor)$ obeys a Gaussian distribution. That is, for $n$ chosen uniformly at random from $\{1,\ldots,N\}$, the random variable $(\omega(\lfloor\alpha n^2\rfloor) - \log\log N)/\sqrt{\log\log N}$ converges in distribution to the standard Gaussian as $N$ tends to infinity. Additionally, we fully describe the distribution of $\omega(\lfloor\alpha n^2\rfloor)$ for a positive rational $\alpha$. 
\end{abstract}
\tableofcontents
\section{Introduction}
The distribution of prime factors of integers is a classical topic in Probabilistic Number Theory. For a natural number $n,$ let $\omega(n) = \sum_{p\mid n} 1$ denote the number of distinct prime divisors of $n.$ By convention, we write $\omega(0) = 0$ and $\omega(-n)=\omega(n)$. Hardy and Ramanujan \cite{HardyRamanujan1920} have established that $\omega(n)$ has normal order $\log\log n.$ That is, for any $\varepsilon > 0$, there exists an $x_\varepsilon$ such that if $x \ge x_\varepsilon$ is sufficiently large, we have 
$$(1-\varepsilon)\log\log x \le \omega(n) \le (1+\varepsilon)\log\log x$$ for all but at most $\varepsilon x$ integers $n\le x$. Erd\H{o}s and Kac \cite{ErdosKac1940} later improved this result by showing that the fluctuations around this normal order are asymptotically Gaussian. More precisely, for every $\tau\in\mathbb R$,
\[
\lim_{N\to\infty}\frac1N
\# \left\{n\le N:
\frac{\omega(n)-\log\log N}{\sqrt{\log\log N}}
\le\tau\right\}
=\Phi(\tau),
\]
where $\Phi$ is the standard Gaussian distribution function
$$\Phi(\tau) = \frac{1}{\sqrt{2\pi}}\int_{-\infty}^{\tau} e^{-t^2/2} \, dt.$$
For an exposition of the Erd\H{o}s-Kac theorem using moments and a sieve setting, one may see Granville and Soundararajan \cite{Granville01}.

A natural extension is to ask if a similar result holds along polynomial sequences instead of all natural numbers. To this end, Minoru Tanaka established in \cite{Tanaka55} an analogue of the Erd\H{o}s-Kac theorem, in which the number of distinct irreducible factors of an integer polynomial determines the mean and variance of the resulting Gaussian law. 

Another natural way of extending the Erd\H{o}s-Kac theorem is by considering  the integer part of real polynomials. Banks and Shparlinski \cite{BANKS2007413} investigated the linear case, that is, the distribution of prime divisors in Beatty sequences of the form $\lfloor\alpha n+\beta\rfloor, $ where $\alpha>0$ is irrational and $\beta$ is a real number. More recently, Fredy Yip \cite{Yip26} obtained quantitative bounds on Erd\H{o}s-Kac laws for Beatty sequences. It has been conjectured that the same Gaussian law holds for $\lfloor f(n)\rfloor$, where $f$ is a polynomial with positive leading coefficient and at least one irrational non-constant coefficient.

Motivated by these developments, we consider the quadratic polynomial
\[
f(n)=\alpha n^2
\]
for $\alpha>0$ and study the asymptotic distribution of $\omega(\lfloor f(n)\rfloor)$. We establish an Erd\H{o}s-Kac law for every irrational $\alpha>0$, without imposing any Diophantine conditions on $\alpha$, and determine for which rationals $\alpha$ such a law holds.

More precisely, we will show the following.
\begin{theorem}\label{thm:main theorem}
    Let $\alpha > 0$ be an irrational number. We have
    $$\lim_{N\to\infty} \frac{1}{N} \#\left\{n \le N : \frac{\omega(\lfloor \alpha n^2 \rfloor)-\log\log N}{\sqrt{\log\log N}} < \tau\right\} = \Phi(\tau),$$
    for any $\tau\in\R,$ where $\Phi$ is the normal Gaussian distribution function.
\end{theorem}

In fact, we will deduce this from an estimate of the corresponding moments. Recall that a standard normal random variable $Z \sim \mathcal{N}(0,1)$ has moments
\[\mathbb E(Z^{\ell}) = \begin{cases}
(\ell-1)!!, & \text{for}\;\ell \;\text{even}, \\
0, & \text{for}\;\ell \;\text{odd}.
\end{cases}\]

\begin{theorem}\label{thm:moment theorem}   
    Let $\alpha > 0$ be irrational. For all $\ell \in \N$ and for every sufficiently large $N \in \N$, there is a set $\mathcal{C}_N \subseteq [N]$ of size $N - o(N)$, independent of $\ell$, such that
    \[\frac{1}{\vert \mathcal{C}_N \vert} \sum_{n \in \mathcal{C}_N} \bigg( \frac{\omega(\lfloor \alpha n^2 \rfloor) - \log \log N}{\sqrt{\log \log N}} \bigg)^{\ell} = \begin{cases}
    (\ell-1)!!, & \text{for}\;\ell \;\text{even}, \\
    0, & \text{for}\;\ell \;\text{odd,}
    \end{cases} + O_{\ell} \bigg( \frac{\log \log \log N}{\sqrt{\log \log N}} \bigg).\]    
\end{theorem}

\begin{proof}[Proof that \Cref{thm:moment theorem} implies \Cref{thm:main theorem}]
    Note that \Cref{thm:main theorem} is equivalent to
    \[\lim_{N\to\infty} \frac{1}{\vert \mathcal{C}_N \vert} \#\left\{n \in \mathcal{C}_N : \frac{\omega(\lfloor \alpha n^2 \rfloor)-\log\log N}{\sqrt{\log\log N}} < \tau\right\} = \Phi(\tau)\]
    for every $\tau \in \R$. For $n$ chosen uniformly at random from $\mathcal{C}_N$, consider the random variable
    \[X_N = \frac{\omega(\lfloor \alpha n^2 \rfloor)-\log\log N}{\sqrt{\log \log N}}.\]
    Clearly, for all sufficiently large $N \in \N$, every moment of $X_N$ exists, and \Cref{thm:moment theorem} implies
    \[\lim_{N \to \infty} \mathbb{E}(X_N^{\ell}) = \mathbb{E}(Z^{\ell}), \qquad Z \sim \mathcal{N}(0,1),\]
    for every $\ell \in \N$. Since the standard normal distribution is well-known to be determined by its moments (\cite[Example 30.1]{BillingsleyProbMeasure}), the method of moments (\cite[Theorem 30.2]{BillingsleyProbMeasure}) thus implies convergence in distribution, proving \Cref{thm:main theorem}.
\end{proof}

By the original Erd\H{o}s-Kac theorem, the conclusion of \Cref{thm:main theorem} also holds for $\alpha \in \N$, as in this case $\omega(\lfloor \alpha n^2 \rfloor) = \omega(n) + O_{\alpha}(1)$. For rational $\alpha$, however, the situation is more subtle, as our third theorem shows. Let $\sqf(n)$ denote the square-free part of $n$.

\begin{theorem}\label{thm:rational theorem}
    Let $\alpha = \frac{a}{q} \in \Q\setminus \N$ be a positive rational number with $(a,q)=1$. The number of distinct prime factors of $\lfloor \alpha n^2\rfloor$ obeys a Gaussian law if and only if $\sqf(a) > q$. 
    In this case, we have
    $$\lim_{N\to\infty} \frac{1}{N} \#\left\{n \le N : \frac{\omega(\lfloor \alpha n^2 \rfloor)-\log\log N}{\sqrt{\log\log N}} < \tau\right\} = \Phi(\tau),$$
    for any $\tau\in\R.$
\end{theorem}

We fully solve the case of rational $\alpha$ by also proving the following.

\begin{corollary}\label{cor: rational corollary}
    Let $\alpha = \frac{a}{q} \in \Q\setminus \N$ be a positive rational number with $(a,q)=1$. The actual distribution of prime factors for $\lfloor\alpha n^2\rfloor$ is that for any $\tau\in\R$, we have
    $$\lim_{N\to\infty} \frac{1}{N} \#\left\{n \le N : \frac{\omega(\lfloor \alpha n^2 \rfloor)-(1 + 1_{S}(n))\log\log N}{\sqrt{(1+1_{S}(n))\log\log N}} < \tau\right\} = \Phi(\tau),$$
    where $S$ is a set consisting of 
    $$\sum_{k\le x}\#\{0\le r<q:r^2\equiv k^2\bmod q\}$$
    many congruence classes mod $q$, where 
    $$x = \left\lfloor\sqrt{\frac{q-1}{\sqf(a)}}\right\rfloor.$$
\end{corollary}

\subsection{Proof Outline}

In order to prove \Cref{thm:moment theorem}, we consider two cases, based on the Diophantine properties of $\alpha$. 

Given a real number $\alpha$, an approximation parameter $0 < \delta < 1/100$, and $N\in\N$, we say that $\alpha$ can be ``well-approximated'' at scale $N$ by a rational if there exists $(a,q) = 1$ with 
$$q \le N^\delta \qquad \text{ and } \qquad \left|\alpha - \frac{a}{q}\right| < N^{\delta-2}.$$
First, if $\alpha$ cannot be ``well-approximated'' by rationals, we say that $\alpha$ is in the minor arcs, and in this case we have a rather standard proof that $\omega(\lfloor\alpha n^2\rfloor)$ obeys a Gaussian law. Indeed, we use the sieve setting of Granville and Soundararajan in \cite{Granville01} combined with the inequalities of Erd\H os--Turán and Weyl to bound the resulting exponential sums. 

On the other hand, if $\alpha$ can be ``well-approximated'' by rationals, we say that $\alpha$ is in the major arcs.
In this case, we proceed as follows. First, we rewrite
$$\frac{1}{\mathcal C_N}\sum_{n\in \mathcal C_N} \left(\omega(\lfloor\alpha n^2\rfloor)-\log\log N \right)^\ell$$
using $\alpha = \frac{a}{q}+\beta$, where $\beta$ is small, to get a sum over irreducible polynomials $f_r$ (with $0\le r < q$).  

To show that $\omega(f_r(m))$ obeys a Gaussian law, we begin by showing that $\omega_z(f_r(m))$ obeys such a law, where $\omega_z$ is a truncation over primes in the range $p \le z$ with $p \nmid 2aq$. To show this, we extract the moments of $\omega_z(f_r(m))$ through the generating function $\exp(t\omega_z(f_r(m)))$. Further, we count the solutions to $p \mid f_r(m)$ using that $f_r(m) \equiv 0 \bmod p$ has $1+\Legendre{\Delta_r}{p}$ solutions $\bmod$ $p$, where $\Delta_r$ is the discriminant of $f_r$.

Completing the proof that $\omega_z(f_r(m))$ obeys a Gaussian law involves showing that the average over all relevant discriminants of 
\begin{equation}\label{eq:Proof plan 1}\sum_{p\le z} \frac{1}{p}\Legendre{\Delta_r}{p}\end{equation}
is small. 
We present a proof that utilizes the Quasi-Riemann Hypothesis for $\Re(s) > 11/12$, which allows us to have very good control over a character sum over primes closely related to (\ref{eq:Proof plan 1}). 

Finally, we are left to show that switching from $\omega_z(\cdot)$ to $\omega(\cdot)$ does not create any uncontrollable error. For this, we first show that on average over all functions $f_r$, the primes dividing both $f_r(m)$ and $aq$ are few. We then conclude by saying that the big primes $p > z$ also have small contribution. 

Bringing all of this together allows us to get the desired result.

\begin{center}
    \textbf{Notation}
\end{center}

We denote by $\N = \{1,2,3,\dots\}$ the set of natural numbers and $p \in \N$ will always denote a prime number. For the complex exponential, we use the convention $e(x) = e^{2\pi i x}$, and for $N \in \N$, we write $[N] = \{1,2,\ldots,N\}$. Furthermore, we abuse notation by writing $\log(x) := \log(\max(x,e))$. We also use the Legendre symbol $\Legendre{a}{p}$ defined as
\[\Legendre{a}{p} = \begin{cases} 
   1 & \text{if } a \text{ is a quadratic residue modulo } p, \\
  -1 & \text{if } a \text{ is a quadratic non-residue modulo } p, \\
   0 & \text{if } a \equiv 0 \pmod p.
\end{cases}\]
Regarding asymptotic notation, we write $X = O(Y)$ or $X \ll Y$ if there is a constant $C$ such that $|X| \le CY$. 
\begin{remark}
    We often use that $\sum_{p_1\cdots p_k \mid n} \frac{1}{p_1\cdots p_k} \ll_k (\log\log\log n)^k.$ This can be deduced from Mertens' theorem combined with the fact that $\omega(n) \ll \log n$ and that $1/p$ is decreasing. In particular, for $k=1$, we have 
    $$\sum_{p \mid n} \frac{1}{p} \ll \log\log\log n.$$
\end{remark}

\begin{center}
    \textbf{Acknowledgments}
\end{center}
We would like to thank Florian K.~Richter for his guidance and support throughout this project. We are also grateful to Nadia Kaiser for her assistance with the organization and logistics of the Young Researchers in Mathematics Program 2026. Finally, we thank the Bernoulli Center for Fundamental Studies and everyone involved in making the program possible for providing us with this opportunity and fostering a stimulating and collaborative research environment.
\begin{center}
    \textbf{AI use disclosure}
\end{center}
In preparing this paper, the authors used AI solely for literature searches, proofreading grammar and punctuation, and checking mathematical accuracy. AI was not used to generate mathematical ideas.

\bigskip
\section{Minor Arcs}\label{sec: Minor arcs}
In what follows, we let $\alpha >0$ be irrational, $0 < \delta < 1/100$, and take $N \in \N$ to be large. Recall that $\alpha$ is in the minor arcs if for every $(a,q)=1$ we have 
$$\left|\alpha - \frac{a}{q}\right| < N^{\delta-2} \implies q > N^\delta.$$
The goal of this section is to prove the following proposition, which establishes \Cref{thm:moment theorem} in this case with $\mathcal{C}_N = [N]$. 
\begin{proposition}[Erd\H os-Kac for Minor Arcs]\label{prop:Minor Arcs}
Let $\ell\in \N$ be fixed. Then 
\[
\frac{1}{N}\sum_{n\le N}\left(\frac{\omega(\lfloor\alpha n^2\rfloor) - \log\log N}{\sqrt{\log\log N}}\right)^\ell=\begin{cases}
(\ell-1)!!, & \text{for}\;\ell \;\text{even},\\
0, & \text{for}\;\ell \;\text{odd},
\end{cases}
+O_\ell\left(\frac{1}{\sqrt{\log\log N}}\right).\]
\end{proposition}
In the proof, we first estimate the moments of a truncated version of $\omega$ before removing this truncation and showing that we obtain a suitably small error. To do the first step, we follow the method that Granville and Soundararajan use in \cite{Granville01}, so we need to reintroduce some of their notation. Let $\mathcal{A}=\{a_1,\dots, a_N\}$ be a multiset of $N$ natural numbers and $\mathcal{P}$ be a set of primes. Define $\mathcal{A}_d=\#\{n\le N\colon d\mid a_n\}$, and suppose that there exists some real-valued, nonnegative multiplicative function $h(d)$, with $0\le h(d)\le d$ for square-free $d$, such that
\[
\mathcal{A}_d=\frac{h(d)}{d}N+r_d,
\]
for square-free $d$, where $r_d$ is expected to be small. Finally, let $\omega_{\mathcal{P}},\mu_{\mathcal{P}}$ and $\sigma_{\mathcal{P}}$ be defined as
\[
\omega_{\mathcal{P}}(a)=\sum_{p\in \mathcal{P}}
1_{p\mid a},\quad \mu_{\mathcal{P}}=\sum_{p\in \mathcal{P}}\frac{h(p)}{p},\quad \text{and}\quad\sigma_{\mathcal{P}}^2=\sum_{p\in \mathcal{P}}\frac{h(p)}{p} \left(1-\frac{h(p)}{p}\right),
\]
for $a\in \mathcal{A}$. The following result from \cite{Granville01} will assist us in proving \Cref{prop:Minor Arcs}.
\begin{lemma}[{\cite[Proposition~3]{Granville01}}]\label{lem:Gran-Sound}
Uniformly for all natural numbers $k\le \sigma_{\mathcal{P}}^{2/3}$, we have
\[
\sum_{a\in \mathcal{A}}\left(\omega_{\mathcal{P}}(a)-\mu_{\mathcal{P}}\right)^k = C_k N\sigma_{\mathcal{P}}^k \left(1 + O\left(\frac{k^3}{\sigma_{\mathcal{P}}^2}\right)\right) + O\left(\mu_{\mathcal{P}}^k \sum_{d\in D_k(\mathcal{P})}|r_d|\right),
\]
if $k$ is even, and 
\[
\sum_{a\in \mathcal{A}}\left(\omega_{\mathcal{P}}(a)-\mu_{\mathcal{P}}\right)^k \ll C_k N\sigma_{\mathcal{P}}^k\frac{k^{\frac 32}}{\sigma_{\mathcal{P}}} +\mu_{\mathcal{P}}^k \sum_{d\in D_k(\mathcal{P})}|r_d|,
\]
if $k$ is odd, where 
\[
C_k=\frac{\Gamma(k+1)}{2^{k/2}\Gamma(k/2+1)},
\]
and $D_k(\mathcal{P})$ denotes the set of square-free integers that are the product of at most $k$ primes, all from the set $\mathcal{P}$.
\end{lemma}
Fix an $\ell \in \N$. We shall apply \Cref{lem:Gran-Sound} for a specific choice of $\mathcal{A}$ and $\mathcal{P}$. To begin with, let $n_0$ be the largest integer such that $\lfloor \alpha n^2\rfloor=0$. We may now define $\mathcal{A}(N)=\{\lfloor \alpha (n+n_0)^2\rfloor\colon n\le N\}$ and $\mathcal{P}_\ell=\{p\colon p\le z_\ell\}$, with $z_\ell$ to be chosen below. We define $\mathcal{A}(N)$ in this way, to avoid any potential issues regarding $\omega_{\mathcal{P}_\ell}(0)$. Then, as in \cite{Granville01}, we define $\mathcal{A}_d(N)=\#\{n\le N\colon d\mid \lfloor\alpha (n+n_0)^2\rfloor\}$.
Observe now that
\[
\mathcal{A}_d(N)=\#\left\{n\le N\colon  \left\{\frac{\alpha (n+n_0)^2}{d}\right\}\in\left[0,\frac{1}{d}\right)\right\}.
\]
Since, for any $d\in \N$, $\alpha/d$ is irrational, Weyl's Equidistribution Theorem yields
\[
\mathcal{A}_d(N)=\frac{1}{d}N+r_d(N),
\]
where $r_d(N)=o_d(N)$.

In order to apply \Cref{lem:Gran-Sound} in our setting, we need an estimate for $\sum_{d\in D_k(\mathcal{P}_\ell)}|r_d(N)|$.
\begin{lemma}\label{lem:Small $r_d$ in the minor arcs}
    Let $\ell\in \N$ and take $z_\ell=N^{\delta/8\ell}$. Then
    $$\sum_{d\in D_k(\mathcal{P}_\ell)}|r_d(N)| = O(N^{1-\delta/8}),$$
    uniformly for $k\le \ell$.
\end{lemma}
The proof of this lemma relies on applying the Erd\H os-Turán inequality to reduce the problem to estimating an exponential sum, for which we use the Weyl inequality. We recall both below for convenience.
\begin{lemma}[{\cite[Chapter 2, Theorem 2.5]{uniformdistributionofsequences}}]\label{lem:Erdos-Turan}
Let $s_1,s_2,s_3,\ldots\in \R$ be a sequence. There exists an absolute constant $C>0$ such that for any $M,N\in \N$
\[
\sup_{0\le a\le b\le 1}\left|\frac{1}{N}\#\{1\le n\le N\colon a\le s_n\mod 1\le b\}-(b-a)\right|\le C\left(\frac{1}{M}+\sum_{m\le M}\frac{1}{m}\left|\frac{1}{N}\sum_{n\le N}e(ms_n)\right|\right).
\]
\end{lemma}
\begin{lemma}[{\cite[Theorem 4.3]{nathanson1996additive}}]
If $M,N,a,q$ are integers with $a,q$ coprime, $q>0$, and $f$ is a real polynomial of degree $k$ with leading coefficient $c$ that satisfies $|c-a/q|\le tq^{-2}$ for some $t\ge 1$, then for any $\varepsilon>0$, we have
\[
\sum_{M\le x\le M+N}e(f(x))=O\left(N^{1+\varepsilon}\left(\frac{t}{q}+\frac{1}{N}+\frac{t}{N^{k-1}}+\frac{q}{N^k}\right)^{2^{1-k}}\right),
\]
as $N\to \infty$.
\end{lemma}
\begin{proof}[Proof of \Cref{lem:Small $r_d$ in the minor arcs}]
Let $d\in D_k(\mathcal{P}_\ell)$. Since 
\[
\frac{1}{N}r_d(N)=\frac{1}{N}\#\left\{n\le N\colon  0\le\left\{\frac{\alpha (n+n_0)^2}{d}\right\}<\frac{1}{d}\right\}-\frac{1}{d},
\]
we can apply the Erd\H os-Turán inequality to the sequence $\{(n+n_0)^2\alpha/d\}$, which yields
\[|r_d(N)| \le C\left(\frac{N}{M} + \sum_{m\le M}\frac{1}{m}\left|\sum_{n\le N}e\left(\frac{m(n+n_0)^2\alpha}{d}\right)\right|\right),
\]
for any $M\in \N$. We take $M = \lfloor N^{\delta/4}\rfloor$, and we want to estimate the exponential sum
\[
\sum_{n\le N} e\left(\frac{m(n+n_0)^2\alpha}{d}\right).
\]
Applying Dirichlet's approximation theorem to $m\alpha/d$, we get that there exist $a',q' \in \N$ with $q'\le N^{2-3\delta/4}$ and $(a',q')=1$ such that 
\[
\left|\frac{m\alpha}{d}-\frac{a'}{q'}\right|<\frac{1}{q'N^{2-3\delta/4}},
\]
which yields
\[
\left|\alpha-\frac{a'd}{q'm}\right|<\frac{d}{q'mN^{2-3\delta/4}}.
\]
Since $d\le z_\ell^k \le N^{\delta/8}$ and $q'm\ge1$, we obtain
\[    
\left|\alpha-\frac{a'd}{q'm}\right|<\frac{1}{N^{2-\delta}}.
\]
But since $\alpha$ is in the minor arcs, we necessarily have $q'm\ge \frac{q'm}{(q'm,a'd)}> N^\delta$, so $q'\ge N^{3\delta/4}$.
Applying the quadratic Weyl's inequality for the polynomial in $n$, $m(n+n_0)^2\alpha/d$, which has leading coefficient $m\alpha/d$, we obtain
\[
\left|\sum_{n\le N} e\left(\frac{m(n+n_0)^2\alpha}{d}\right)\right| \ll_{\varepsilon} N^{1+\varepsilon}\left(\frac{1}{q'} + \frac{1}{N} + \frac{q'}{N^2}\right)^{1/2} \ll_{\varepsilon} N^{1+\varepsilon-3\delta/8},
\]
for every $\varepsilon>0$. Substituting this bound into the Erd\H os-Turan inequality, we have 
\[
|r_d(N)| \ll_{\varepsilon} \left(N^{1-\delta/4} + \log (N) N^{1+\varepsilon-3\delta/8}\right) \ll N^{1-\delta/4},
\]
using $\log N \ll_\varepsilon N^{\varepsilon}$ and taking $\varepsilon = \delta / 16$. This implies that 
\[
\sum_{d\in D_k(\mathcal{P}_\ell)} |r_d(N)| \ll |D_k(\mathcal{P}_\ell)|\cdot N^{1-\delta/4} \ll N^{k\delta/8\ell}\cdot N^{1-\delta/4} \ll N^{1-\delta/8},
\]
which concludes our proof. 
\end{proof}
 
Applying \Cref{lem:Gran-Sound} for our specific choice of $\mathcal{A}$ and $\mathcal{P}_\ell$ and using \Cref{lem:Small $r_d$ in the minor arcs}, we get the following result.
\begin{corollary}\label{cor:GS for minor arcs}
 Let $\ell\in \N$. Then uniformly for $k \le \min\{\ell,\sigma_{\mathcal{P}_\ell}^{\frac{2}{3}}\}$
\[
\sum_{n\le N}\left(\omega_{\mathcal{P}_\ell}(\lfloor\alpha (n+n_0)^2\rfloor) - \mu_{\mathcal{P}_\ell}\right)^k = C_k N\sigma_{\mathcal{P}_\ell}^k \left(1 + O\left(\frac{k^3}{\sigma_{\mathcal{P}_\ell}^2}\right)\right) + O_{\ell}\left(\mu_{\mathcal{P}_\ell}^k N^{1-\delta/8}\right),
\]
if $k$ is even, and 
\[
\sum_{n\le N}\left(\omega_{\mathcal{P}_\ell}(\lfloor\alpha (n+n_0)^2\rfloor) - \mu_{\mathcal{P}_\ell}\right)^k \ll_\ell C_kN\sigma_{\mathcal{P}_\ell}^k\frac{k^{\frac 32}}{\sigma_{\mathcal{P}_\ell}} + \mu_{\mathcal{P}_\ell}^k N^{1-\delta/8},
\]
if $k$ is odd, 
where 
\[
C_k=\frac{\Gamma(k+1)}{2^{k/2}\Gamma(k/2+1)}.
\]
\end{corollary}
We can now proceed to prove \Cref{prop:Minor Arcs}.
\begin{proof}[Proof of \Cref{prop:Minor Arcs}]
    What remains is to remove the truncation of $\omega$ and transition to the usual moments. Fix an $\ell\in \N$ and $N$ large enough so that $\ell\le\sigma_{\mathcal{P}_\ell}^{\frac{2}{3}}$. If we define $U_N(n)=\omega(\lfloor\alpha (n+n_0)^2\rfloor)-\log\log N$ and $V_N(n)=\omega_{\mathcal{P}_\ell}(\lfloor\alpha (n+n_0)^2\rfloor)-\mu_{\mathcal{P}_\ell}$, then we get
\[
U_N=V_N+(\omega(\lfloor\alpha (n+n_0)^2\rfloor)-\omega_{\mathcal{P}_\ell}(\lfloor\alpha (n+n_0)^2\rfloor))+(\mu_{\mathcal{P}_\ell}-\log\log N).
\]

Using the binomial theorem, we obtain
\begin{align*}
    \frac{1}{N}\sum_{n\le N}U_N^\ell&= \frac{1}{N}\sum_{n\le N}\sum_{k=0}^{\ell}\binom{\ell}{k}V_N^{\ell-k}(U_N-V_N)^k\\
    &= \frac{1}{N}\sum_{n\le N}V_N^\ell+\sum_{k=1}^{\ell}\binom{\ell}{k} \frac{1}{N}\sum_{n\le N}V_N^{\ell-k}\big(U_N-V_N\big)^k.
\end{align*}
We first establish a uniform bound for the difference $U_N-V_N$. First,
\[
|\omega(\lfloor\alpha (n+n_0)^2\rfloor)-\omega_{\mathcal{P}_\ell}(\lfloor\alpha (n+n_0)^2\rfloor)|\le \frac{\log \lfloor\alpha (n+n_0)^2\rfloor}{\log z_\ell}=O_{\ell}(1).
\]
Moreover, by Mertens' theorem,
\[
\mu_{\mathcal{P}_\ell}-\log\log N=O_{\ell}(1), \quad \text{and}\quad \sigma_{\mathcal{P}_\ell}^2-\log\log N=O_{\ell}(1),
\]
and therefore $U_N-V_N=O_\ell(1)$. Furthermore, we observe that 
\[ \frac{1}{N}\sum_{n\le N}|V_N|^j\ll_{\ell}\sigma_{\mathcal{P}_\ell}^{j}
\]
for every $0\le j\le \ell-1$. Indeed, for $j=0$, it is trivial, so let $1\le j\le \ell-1$ be even. Then \[
\frac{1}{N}\sum_{n\le N}|V_N|^j= \frac{1}{N}\sum_{n\le N}V_N^j,
\]
and so, by \Cref{cor:GS for minor arcs},
\[
 \frac{1}{N}\sum_{n\le N}V_N^j=\frac{1}{N}\sum_{n\le N}(\omega_{\mathcal{P}_\ell}(\lfloor\alpha (n+n_0)^2\rfloor)-\mu_{\mathcal{P}_\ell})^j\ll_\ell\sigma_{\mathcal{P}_\ell}^{j}.
\]
For odd $j$, by the Cauchy-Schwarz inequality, we obtain
\[
 \frac{1}{N}\sum_{n\le N}|V_N|^j\le  \left(\frac{1}{N}\sum_{n\le N}V_N^{j-1}\right)^{1/2}\left(\frac{1}{N}\sum_{n\le N}V_N^{j+1}\right)^{1/2}\ll_\ell\sigma_{\mathcal{P}_\ell}^{j},
\]
because $j-1,j+1\le \ell$. Consequently, for $1\le k\le \ell$, we deduce that
\begin{align*}
    \bigg|\frac{1}{N}\sum_{n\le N}V_N^{\ell-k}\left(U_N-V_N\right)^k\bigg|&\le \sup_{n\le N}|U_N-V_N|^k\frac{1}{N}\sum_{n\le N}|V_N|^{\ell-k}\ll_\ell\sigma_{\mathcal{P}_\ell}^{\ell-k}.
\end{align*}
Hence, we have
\[
\sum_{k=1}^{\ell}\binom{\ell}{k}\frac{1}{N}\sum_{n\le N}V_N^{\ell-k}(U_N-V_N)^k\ll_\ell \sum_{k=1}^{\ell}\binom{\ell}{k}\sigma_{\mathcal{P}_\ell}^{\ell-k}\ll_\ell \sigma_{\mathcal{P}_\ell}^{\ell-1}\ll_\ell(\log\log N)^{(\ell-1)/2},
\]

and therefore
\begin{align*}
    \frac{1}{N}\sum_{n\le N}\left(\frac{\omega(\lfloor\alpha (n+n_0)^2\rfloor) - \log\log N}{\sqrt{\log\log N}}\right)^\ell=\frac{1}{N}\sum_{n\le N}\left(\frac{\omega_{\mathcal{P}_\ell}(\lfloor \alpha (n+n_0)^2\rfloor)-\mu_{\mathcal{P}_\ell}}{\sigma_{\mathcal{P}_\ell}}\right)^\ell\!\Bigg(1& +O_\ell\left(\frac{1}{\sqrt{\log\log N}}\right)\!\Bigg)\\
    &+O_\ell\left(\frac{1}{\sqrt{\log\log N}}\right),
\end{align*}
where, when normalizing by $\sqrt{\log\log N}$, we use that 
\[
\frac{\sigma_{\mathcal{P}_\ell}}{\sqrt{\log\log N}}=1+O_\ell\left(\frac{1}{\sqrt{\log\log N}}\right).
\]
To conclude, we first apply \Cref{cor:GS for minor arcs}, and use that, for $\ell$ even, $C_\ell=(\ell-1)!!$, and that
\[
\frac{1}{\sigma_{\mathcal{P}_\ell}^2},\frac{1}{\sigma_{\mathcal{P}_\ell}},\left(\frac{\mu_{\mathcal{P}_\ell}}{\sigma_{\mathcal{P}_\ell}}\right)^\ell\frac{1}{N^{\delta/8}}\ll_\ell\left(\frac{1}{\sqrt{\log\log N}}\right).
\]
Finally, if we observe that in the above we have just shifted our original sequence by $n_0$, we get \Cref{prop:Minor Arcs}.
\end{proof}
\section{Major arcs}
As before, we let $\alpha > 0$ be irrational, $0 < \delta < 1/100$, and take $N \in \N$ to be large. In the previous section, we dealt with the case of $\alpha$ being in the minor arcs, so for what follows we assume that $\alpha$ lies in the major arcs. Therefore, we let $a,q\in \N$, with $(a,q)=1$, set such that
$$q \le N^\delta \qquad \text{ and } \qquad \left|\alpha - \frac{a}{q}\right| < N^{\delta-2}.$$
Here, it is to be understood that $\alpha$ can be ``well approximated'' by a rational number. This will be particularly useful to convert the following into integer polynomials, which are well understood, notably through Minoru Tanaka's proof in \cite{Tanaka55}. The goal will thus be to prove the following  version of \Cref{thm:moment theorem} for the major arcs.

\begin{proposition}[Erd\H os--Kac for Major Arcs]\label{prop:Major Arcs}
Let $\alpha$ be in the major arcs and $\ell\in \N$ be fixed. Then
$$\frac{1}{\vert \mathcal{C}_N \vert}\sum_{n \in \mathcal{C}_N}\left(\frac{\omega(\lfloor\alpha n^2\rfloor) - \log\log N}{\sqrt{\log\log N}}\right)^\ell=\begin{cases}
(\ell-1)!! ,& \text{ for }\ell \;\text{ even },\\
0, & \text{ for }\ell \text{ odd },
\end{cases}
+O_\ell \left( \frac{\log \log \log N}{\sqrt{\log \log N}} \right).$$
\end{proposition}
\subsection{From \texorpdfstring{$\lfloor\alpha n^2\rfloor$}{floor(alpha n squared)} to Polynomials}\label{subsec: alpha to polynomials}
To prove the above statement, we start by converting the sum over $\lfloor\alpha n^2\rfloor$ to a sum over polynomials. This is the goal of the following proposition. At a later stage, we will ensure the polynomials are irreducible. 

\begin{proposition}\label{prop:from floor to polynomial}
    Let $K = \left\lfloor\frac{\sqrt N}{q}\right\rfloor$. For every $\ell \in \N$, there exist intervals $I_{n,r}^{-},I_{n,r}^{+} \subseteq \{0,1,\ldots, K\}$ of size at least $N^{1/4}$ such that we bound
    \begin{multline*}
        \frac{qK}{N}\sum_{0\le n \le \frac{N}{qK}-1}\min_{\eta\in\{0,\pm1\}^q}\left\{\frac{1}{q} \sum_{0 \leq r < q}\frac{1}{|I_{n,r}^-|}\sum_{m \in I_{n,r}^-}\left(\omega(f_{n,r,\eta_r}(m)) - \log\log N\right)^\ell\right\}+o_\ell(1)\\
        \le \frac{1}{N}\sum_{n\in [N]} \left(\omega(\lfloor\alpha n^2\rfloor)-\log\log N \right)^\ell \le \\
        \frac{qK}{N}\sum_{0\le n \le \frac{N}{qK}-1}\max_{\eta\in\{0,\pm1\}^q}\left\{\frac{1}{q} \sum_{0 \leq r < q}\frac{1}{|I_{n,r}^+|}\sum_{m \in I_{n,r}^+}\left(\omega(f_{n,r,\eta_r}(m)) - \log\log N\right)^\ell\right\}+o_\ell(1),
    \end{multline*}
    with 
    $$f_{n,r,\eta_r}(m) = aqm^2 +(2nKaq + 2ar)m + c_{n,r}-\eta_r$$
    such that 
    $$c_{n,r} = \lfloor \alpha(nqK)^2 + \frac{ar^2}{q}\rfloor + 2nKar \qquad \text{ and } \qquad \eta = (\eta_0, \dots, \eta_{q-1}) \in \{0,\pm1\}^q.$$
 \end{proposition}
 \begin{proof}
     First, we split the sum over ranges of size $qK$, so roughly of size $\sqrt{N}$. We obtain
     \begin{multline*}
         \frac{1}{N}\sum_{n\le N} (\omega(\lfloor\alpha n^2\rfloor)-\log\log N)^\ell = \frac{1}{N}\sum_{n = 0}^{\lfloor\frac{N}{qK}-1\rfloor}\sum_{k \le qK} (\omega(\lfloor\alpha (nqK + k)^2\rfloor)-\log\log N)^\ell + O\left(\frac{qK(\log N)^\ell}{N}\right) \\
         = \frac{qK}{N}\sum_{n = 0}^{\lfloor\frac{N}{qK}-1\rfloor}\frac{1}{qK}\sum_{k \le qK} \left(\omega(2nKak +\lfloor\alpha (nqK)^2 + \frac{ak^2}{q} + \beta(2nqKk + k^2)\rfloor)-\log\log N\right)^\ell + o_\ell(1),
     \end{multline*}
     where we write $\beta = \alpha - \frac{a}{q}$. 
     Using that $|\beta|(2nqKk + k^2) < 4N^{\delta-2}N^{3/2} < \frac{1}{q}$, we have that $\beta(2nqKk + k^2)$ changes the value inside $\omega(\cdot)$ by either 0, 1, or $-1$. Let us call this shift $\eta_{n,k}$. 
     Additionally, we write $k = mq + r$ for $0\le r<q$, which gives us 
    \begin{multline*}
         \frac{qK}{N}\sum_{n = 0}^{\lfloor\frac{N}{qK}-1\rfloor}\frac{1}{q}\sum_{0\le r < q}\frac{1}{K}\cdot\\
         \sum_{0\le m <K} \left(\omega(2nKa(qm+r) + aqm^2 + 2arm + \lfloor \alpha(nqK)^2 + \frac{ar^2}{q} - \eta_{n,k} \rfloor)-\log\log N\right)^\ell + o_\ell(1),
    \end{multline*}
    which equals 
    $$\frac{qK}{N}\sum_{n = 0}^{\lfloor\frac{N}{qK}-1\rfloor}\frac{1}{q}\sum_{0\le r < q}\frac{1}{K}\sum_{0\le m <K} \left(\omega( aqm^2 +(2nKaq + 2ar)m + c_{n,k})-\log\log N\right)^\ell +o_\ell(1),$$
    where we have 
    $$c_{n,k} = \lfloor \alpha(nqK)^2 + \frac{ar^2}{q}\rfloor + 2nKar -\eta_{n,k},$$ 
    with 
    $$\eta_{n,k} = 
    \begin{cases}
        -1, &\text{ if } \beta >0 \text{ and } 1-\{\alpha(nqK)^2 + ar^2/q\} < \beta(2nqKk+k^2), \\
        1, &\text{ if } \beta <0 \text{ and } \{\alpha(nqK)^2 + ar^2/q\} < |\beta|(2nqKk+k^2), \\
        0, &\text{ otherwise.}
    \end{cases}$$
    
    Notice that within a fixed congruence class modulo $q$ of $k$, $\eta_{n,k}$ is monotone and changes value at most once, which gives rise to two intervals. Let us denote those intervals $I_1 \cup I_2 = \{0,1,\ldots,K-1\}$. If $\eta_{n,k}$ does not change value within a congruence class, we set $I_2 = \emptyset$. We distinguish two cases according to whether $|I_i| < \sqrt[4]{N}$ or not. 
    
    If this is the case, it happens for at most one of them, say $I_2$. In this case we can simply sum and normalize over $I_1$ as this creates an extra error of at most 
    $$\left|\frac{w_1+w_2}{|I_1|+|I_2|}\right| - \left|\frac{w_1}{|I_1|}\right| \le \frac{w_2}{|I_1|+|I_2|}+\frac{w_1}{|I_1|}\cdot\frac{|I_2|}{|I_1|+|I_2|}\ll_\ell2\frac{|I_2|(\log q)^\ell}{K} \ll_\ell o(1),$$
    where $w_i$ represents the sum over $I_i$ and each term within the $\ell$-th power is bounded $\log q$. 
    On the other hand, if all considered intervals have length $|I_i| \ge \sqrt[4]{N}$. We also note that we are able to only sum and normalize over one of the intervals as $$\min\left(\frac{w_1}{|I_1|}, \frac{w_2}{|I_2|}\right) \le \frac{w_1+w_2}{|I_1|+|I_2|} \le \max\left(\frac{w_1}{|I_1|}, \frac{w_2}{|I_2|}\right).$$ 
    We call $I_{n,r}^+$ the interval contributing the most in $\max(\frac{w_1}{|I_1|}, \frac{w_2}{|I_2|})$ and $I_{n,r}^-$ the one contributing the least in $\min(\frac{w_1}{|I_1|}, \frac{w_2}{|I_2|})$. Note that if one of the intervals had length less than $\sqrt[4]{N}$, we have $I_{n,r}^- = I_{n,r}^+ = I_i$, with $|I_i| \ge \sqrt[4]{N}$. 
    
    In particular, we will bound our sum for all possible $\eta_{n,r} \in\{0,\pm1\}$. Setting $c_{n,r} = \lfloor \alpha(nqK)^2 + \frac{ar^2}{q}\rfloor + 2nKar$, we have that our original sum is bounded by 
    \begin{multline*}
        \frac{qK}{N}\sum_{n = 0}^{\lfloor\frac{N}{qK}-1\rfloor}\frac{1}{q}\sum_{0\le r < q}\min_{\eta_r\in\{0,\pm1\}}\left\{\frac{1}{|I_{n,r}^{-}|}\sum_{m \in I_{n,r}^{-}}\left(\omega(f_{n,r,\eta_r}(m)) - \log\log N\right)^\ell\right\} +o_\ell(1) \\
        \le \frac{1}{N}\sum_{n\in[N]} \left(\omega(\lfloor\alpha n^2\rfloor)-\log\log N \right)^\ell \le \\
        \frac{qK}{N}\sum_{n = 0}^{\lfloor\frac{N}{qK}-1\rfloor}\frac{1}{q}\sum_{0\le r < q}\max_{\eta_r\in\{0,\pm1\}}\left\{\frac{1}{|I_{n,r}^{+}|}\sum_{m \in I_{n,r}^{+}}\left(\omega(f_{n,r,\eta_r}(m)) - \log\log N\right)^\ell\right\} +o_\ell(1),
    \end{multline*}
    where $f_{n,r,\eta_r}(m) = aqm^2 +(2nKaq + 2ar)m + c_{n,r}-\eta_r.$
    This is equivalent to taking the max out of $3^q$ different scenarios, yielding
    \begin{equation}\label{eq:result first polynomial}
    \begin{multlined}[0.95\linewidth]
        \frac{qK}{N}\sum_{0\le n \le \frac{N}{qK}-1}\min_{\substack{\eta\in\{0,\pm1\}^q}}\left\{\frac{1}{q} \sum_{0 \leq r < q}\frac{1}{|I_{n,r}^-|}\sum_{m \in I_{n,r}^-}\left(\omega(f_{n,r,\eta_r}(m)) - \log\log N\right)^\ell\right\} +o_\ell(1) \\
        \le \frac{1}{N}\sum_{n\in [N]} \left(\omega(\lfloor\alpha n^2\rfloor)-\log\log N \right)^\ell \le \\
        \frac{qK}{N}\sum_{0\le n \le \frac{N}{qK}-1}\max_{\substack{\eta\in\{0,\pm1\}^q}}\left\{\frac{1}{q} \sum_{0 \leq r < q}\frac{1}{|I_{n,r}^+|}\sum_{m \in I_{n,r}^+}\left(\omega(f_{n,r,\eta_r}(m)) - \log\log N\right)^\ell\right\} +o_\ell(1)
    \end{multlined}
    \end{equation}
    where $\eta_r$ denotes the $r$-th coordinate of $\eta$.
\end{proof}

For what follows, we fix $n$ and $\eta$ and wish to know when $f_{n,r,\eta_r}$ has two distinct irreducible factors. In particular, we denote by $\Delta_{n,r,\eta_r}$ the discriminant of $f_{n,r,\eta_r}$. We have
$$\Delta_{n,r,\eta_r} = (2nKaq+2ar)^2 - 4aq(c_{n,r}-\eta_r) = (2nKaq+2ar)^2 - 4aq(\lfloor \alpha(nqK)^2 + \frac{ar^2}{q}\rfloor + 2nKar - \eta_r).$$ 
Using the property that $\lfloor x \rfloor + \lfloor y\rfloor \le \lfloor x+y\rfloor \le \lfloor x \rfloor + \lfloor y \rfloor + 1$, we can split $\lfloor \alpha(nqK)^2 + \frac{ar^2}{q}\rfloor$, which allows us to control the discriminant with an additional shift $\lambda$ which ranges in $\{0,\pm1,\pm2\}^q$. This yields
\begin{align*}
    \frac{\Delta_{n,r,\eta_r}}{4} &= (nKaq)^2+2(nKaq)(ar)+(ar)^2 - aq\lfloor\alpha(nqK)^2\rfloor - aq\lfloor\frac{ar^2}{q}\rfloor - 2nKa^2qr - aq\lambda_r \\
    &= (nKaq)^2+2(nKaq)(ar)+(ar)^2 - aq\lfloor\alpha(nqK)^2\rfloor - (ar)^2 + aq\{\frac{ar^2}{q}\} - 2nKa^2qr - aq\lambda_r \\
    &= (nKaq)^2 - aq\lfloor \alpha(nqK)^2\rfloor  - aq\lambda_r + a(ar^2 \bmod q),
\end{align*}
where $\lambda_r$ is defined as $\lambda_r = -\eta_r - \left(\lfloor \alpha(nqK)^2\rfloor + \lfloor\frac{ar^2}{q}\rfloor - \lfloor \alpha(nqK)^2 + \frac{ar^2}{q}\rfloor\right)$. We then rewrite $(nKaq)^2 - aq\lfloor \alpha(nqK)^2\rfloor$ using $\alpha = \frac{a}{q}+\beta$ and obtain
\begin{equation}
    \Delta_{n,r,\eta_r} = 4a(ar^2 \bmod q) - 4aq\left(\lfloor\beta(nqK)^2\rfloor + \lambda_r\right).
\end{equation}
We are ready to show that few $f_{n,r,\eta_r}$ have two distinct irreducible factors. 

\begin{lemma}\label{lem: few reducible polynomials}
    Let $X = \{0\le r < q : \Delta_{n,r,\eta_r} \text{ is a square}\}$. We have that for any $\varepsilon > 0$, $|X| \ll_\varepsilon q^{5/6+\varepsilon}$.
\end{lemma}
\begin{proof}
    We will treat each possible value of $\lambda_r \in \{0,\pm1,\pm2\}$ independently. We denote $C_{n,0} = (nKaq)^2-aq\lfloor \alpha(nqK)^2\rfloor$ and $C_{n,\lambda_r} = C_{n,0}-aq\lambda_r$. We then have $\Delta_{n,r,\eta_r}/4 = C_{n,\lambda_r} + a(ar^2 \bmod q)$. We see that $C_{n,\lambda_r} + a(ar^2 \bmod q)$ ranges over at most $q$ values in arithmetic progression (namely $C_{n,\lambda_r} + a\cdot \{0,1,\dots,q-1\}$). By \cite{BombieriSquares}, there are at most $O(q^{2/3}(\log q)^A)$ distinct squares obtainable. The same argument applies for every possible value of $C_{n,\lambda_r}$, giving us a total of at most $O(5q^{2/3}(\log q)^A)$ squares obtainable for a fixed $n$. Let us denote this set by $C_n$. Then, note that as $(a,q)=1$, we can simply count the number of $r$ giving $$(r^2 \bmod q) \in C_n' = \{s : a(as \bmod q)+C_{n,\lambda_r} \in C_n, \text{ for }\lambda_r \in \{0,\pm1,\pm2\}\}.$$
    We now denote by $N_q(s)$ the number of solutions to $r^2 \equiv s \pmod q$. Notice that the number of squares is thus bounded by $\sum_{s\in C_n'} N_q(s)$. Using the Cauchy--Schwarz inequality, we bound this sum by
    $$\sum_{s\in C_n'}N_q(s) \le \left(\sum_{s\in C_n'}1^2\right)^{1/2} \cdot \left(\sum_{s\in C_n'}N_q(s)^2\right)^{1/2} \le |C_n|^{1/2}\left(\sum_{s\in \Z/q\Z}N_q(s)^2\right)^{1/2}.$$

    Observe that $$\sum_{s\in \Z/q\Z}N_q(s)^2=\#\{(x,y) \in (\Z/q\Z)^2 : x^2 \equiv y^2 \pmod q\}.$$
    We can count the pairs in the set using the Chinese Remainder Theorem to count the pairs modulo $p^k \mid\mid q$. Indeed, if $x^2 \equiv y^2 \pmod{p^k}$, we have that $p^k \mid (x-y)(x+y)$. Let us first assume $p>2$. We see that if $p^j \mid(x-y)$ and $p^j \mid (x+y)$, then it divides their sum and difference, so $p^j \mid x$ and $p^j\mid y$. Thus, we write $x = p^j x_1$ and $y = p^j y_1$ with $p^j \mid\mid \gcd(x-y,x+y)$. We split further into two cases: 
    
    \underline{If $k \le 2j$:} we have that $x^2 \equiv y^2 \equiv 0 \pmod{p^k}$. So we have $O(\sqrt {p^k})$ choices for $x$ and $y$, giving us $O(p^k)$ total possible pairs. 

    \underline{If $k > 2j$:} we have $p^{k-2j}\mid (x_1-y_1)(x_1+y_1)$, giving us that either $p^{k-2j} \mid (x_1-y_1)$ or $p^{k-2j} \mid (x_1+y_1)$, hence $x_1 \equiv \pm y_1 \pmod{p^{k-2j}}.$ This means $x$ can be any multiple of $p^{j}$, there are $p^{k-j}$ options. For a fixed $x$, $y$ has 2 options mod $p^{k-2j}$ so $2\frac{p^k}{p^{k-2j}p^j} = 2p^j$ options. This gives $2p^k$ pairs of options. However, as there are up to $k/2$ such $j$, the total number of options is $O(kp^k)$. 

    If $p = 2,$ an extra factor of two appears as $p^k$ may divide $x+y$ and $x-y$ without dividing $x$ and $y$. Thus, the same argument gives a total of $O(2kp^k) = O(kp^k)$ options. Let us call the implied constant $C$. Combining the count for each $p\mid q$ gives us that $$\sum_{s\in \Z/q\Z}N_q(s)^2 = \prod_{p^k\mid\mid q}N(p^k) \le \prod_{p^k\mid\mid q} Ckp^k \ll C^{\omega(q)}q \prod_{p^k\mid\mid q} k \le C^{\omega(q)}qd(q) \ll_\varepsilon q^{1+\varepsilon},$$
    using $\omega(q) = O\left(\frac{\log q}{\log\log q}\right)$ and $d(q) \ll_\varepsilon q^{\varepsilon}$, where $d(q)$ is the number of divisors of $q$. 
    Therefore, we have 
    $$|X|\le\sum_{s\in C_n'}N_q(s) \ll_\varepsilon |C_n'|^{1/2}\left(q^{1+\varepsilon}\right)^{1/2} \le |C_n|^{1/2}\left(q^{1+\varepsilon}\right)^{1/2} \ll_\varepsilon q^{5/6 + \varepsilon}.$$
    Giving us that there are at most $O_\varepsilon(q^{5/6+\varepsilon})$ many values of $r$ for which $\Delta_{n,r,\eta_r}$ is a square. 
\end{proof}
Using \Cref{lem: few reducible polynomials}, we define a set $S \subseteq \{1,2,\ldots,N\}$ of size $o(N)$ such that any integer $m\in [N]$ whose associated polynomial is reducible is included in $S$. Specifically, we define $$S = \bigcup_{n=0}^{\lfloor N/qK\rfloor -1}S_n,$$ where $s \in S_n$ if and only if $s \in \{nqK,nqK+1,\ldots, (n+1)qK-1\}$ and $\Delta_{n,r,x}$ is a perfect square for some $x\in\{0,\pm1\}$, with $r = (s \pmod q)$. We know that $$|S| \le \sum_{n = 0}^{\lfloor\frac{N}{qK}-1\rfloor} |S_n| \ll_\varepsilon \frac{N}{qK}\cdot\frac{qK}{q}\cdot q^{5/6+\varepsilon} = \frac{N}{q^{1/6-\varepsilon}} = o(N),$$
as $q\to\infty$ when $N\to\infty$, using that $\alpha$ is irrational. 

Furthermore, we know, combined with \Cref{eq:result first polynomial} that the sum restricted to $\mathcal C_N = [N]\setminus S$ yields
\begin{multline*}
    \frac{qK}{N}\sum_{0\le n \le \frac{N}{qK}-1}\min_{\substack{\eta\in\{0,\pm1\}^q}}\left\{\frac{1}{q} \sum_{\substack{0 \leq r < q \\ nqK + r \notin S}}\frac{1}{|I_{n,r}^-|}\sum_{m \in I_{n,r}^-}\left(\omega(f_{n,r,\eta_r}(m)) - \log\log N\right)^\ell\right\} +o_\ell(1) \\
    \le\frac{1}{N}\sum_{n\in [N]\setminus S} \left(\omega(\lfloor\alpha n^2\rfloor)-\log\log N \right)^\ell \leq \\
    \frac{qK}{N}\sum_{0\le n \le \frac{N}{qK}-1}\max_{\substack{\eta\in\{0,\pm1\}^q}}\left\{\frac{1}{q} \sum_{\substack{0 \leq r < q \\ nqK + r \notin S}}\frac{1}{|I_{n,r}^+|}\sum_{m \in I_{n,r}^+}\left(\omega(f_{n,r,\eta_r}(m)) - \log\log N\right)^\ell\right\} +o_\ell(1).
\end{multline*}

We note that we are still normalizing by $1/N$, when we would like to normalize by $1/|\mathcal C_N|$. Let $\mathcal{R}_n$ denote the number of residue classes modulo $q$ in the interval $[(n+1)qK-1]\setminus[nqK-1]$ that are not contained in $S_n.$ Since each residue class contains exactly $K$ integers, we have
\[
\mathcal R_n=q-\frac{|S_n|}{K}.
\]
Since the sets $S_n$ are disjoint, we obtain
\begin{align*}
    K\sum_{0\leq n\leq \frac{N}{qK}-1}\mathcal R_n
    &=qK\left\lfloor\frac{N}{qK}\right\rfloor
    -\sum_{0\leq n\leq \frac{N}{qK}-1}|S_n|\\
    &=N-|S|+O(qK).
\end{align*}
And thus,
\[
|\mathcal C_N|=K\sum_{0\leq n \leq \frac{N}{qK}-1}\mathcal R_n+O(qK).
\]
Since $|\mathcal C_N|\sim N$ and $qK\ll\sqrt{N}$, replacing
$|\mathcal C_N|$ by
$K\sum_{0\leq n\leq \frac{N}{qK}-1}\mathcal R_n$
in the normalization introduces an error of
\[
O_\ell \left( \frac{qK}{N}(\log N)^\ell \right)=o_\ell(1).
\]
Then, upon renormalizing, we will get \begin{equation*}
    \begin{multlined}[0.95\linewidth]
        \sum_{0\leq n\leq \frac{N}{qK}-1}
        \frac{\mathcal R_n}
        {\sum_{0\leq j\leq \frac{N}{qK}-1}\mathcal R_j}
        \min_{\substack{\eta\in\{0,\pm1\}^q}}
        \left\{
            \frac{1}{\mathcal R_n}
            \sum_{\substack{0\leq r<q\\ nqK+r\notin S}}
            \frac{1}{|I_{n,r}^-|}
            \sum_{m\in I_{n,r}^-}
            \left(\omega(f_{n,r,\eta_r}(m))
            -\log\log N\right)^\ell
        \right\}+o_\ell(1)\\
        \leq
        \frac{1}{|\mathcal C_N|}
        \sum_{n\in\mathcal C_N}
        \left(\omega(\lfloor\alpha n^2\rfloor)
        -\log\log N\right)^\ell
        \leq\\
        \sum_{0\leq n\leq \frac{N}{qK}-1}
        \frac{\mathcal R_n}
        {\sum_{0\leq j\leq \frac{N}{qK}-1}\mathcal R_j}
        \max_{\substack{\eta\in\{0,\pm1\}^q}}
        \left\{
            \frac{1}{\mathcal R_n}
            \sum_{\substack{0\leq r<q\\ nqK+r\notin S}}
            \frac{1}{|I_{n,r}^+|}
            \sum_{m\in I_{n,r}^+}
            \left(\omega(f_{n,r,\eta_r}(m))
            -\log\log N\right)^\ell
        \right\}+o_\ell(1).
    \end{multlined}
\end{equation*}

Since 
\[
\sum_{0\leq n\leq \frac{N}{qK}-1}\frac{\mathcal R_n}{\sum_{0\leq j\leq \frac{N}{qK}-1}\mathcal R_j}=1
\]
and all the weights are nonnegative, we can now loosen the bounds by taking the minimum and maximum over $n$ as well. This yields:

\begin{equation}\label{eq: double bounds}
    \begin{multlined}[0.95\linewidth]
        \min_{\substack{\eta\in\{0,\pm1\}^q \\ 0 \le n \le \frac{N}{qK}-1}}\left\{\frac{1}{\mathcal R_{n}} \sum_{\substack{0 \leq r < q \\ nqK + r \notin S}}\frac{1}{|I_{n,r}^-|}\sum_{m \in I_{n,r}^-}\left(\omega(f_{n,r,\eta_r}(m)) - \log\log N\right)^\ell\right\} +o_\ell(1) \\
        \le\frac{1}{|\mathcal C_N|}\sum_{n\in \mathcal C_N} \left(\omega(\lfloor\alpha n^2\rfloor)-\log\log N \right)^\ell \leq \\
        \max_{\substack{\eta\in\{0,\pm1\}^q \\ 0 \le n \le \frac{N}{qK}-1}}\left\{\frac{1}{\mathcal R_{n}} \sum_{\substack{0 \leq r < q \\ nqK + r \notin S}}\frac{1}{|I_{n,r}^+|}\sum_{m \in I_{n,r}^+}\left(\omega(f_{n,r,\eta_r}(m)) - \log\log N\right)^\ell\right\} +o_\ell(1).
    \end{multlined}
\end{equation}
In the following subsections, we fix $n$ and $\eta$, and we prove properties within. In this setting, we have
\begin{equation}\label{eq: def f_r}
   f_{r}(m) = aqm^2 +(2nKaq + 2ar)m + c_{n,r}-\eta_r, 
\end{equation}
with discriminant 
\begin{equation}
       \Delta_{r} = 4a\left((ar^2 \bmod q) - q(\lfloor\beta (nqK)^2\rfloor +\lambda_r)\right),
\end{equation}
where we write $K = \lfloor \frac{\sqrt{N}}{q} \rfloor$, $c_{n,r} = \lfloor \alpha(nqK)^2 + \frac{ar^2}{q}\rfloor + 2nKar$, and $$\lambda_r = -\eta_r-\left(\lfloor \alpha(nqK)^2\rfloor + \lfloor\frac{ar^2}{q}\rfloor - \lfloor \alpha(nqK)^2 + \frac{ar^2}{q}\rfloor\right).$$

\subsection{Reduction to Gaussian Moments}
We now show that the moment calculation for the
quadratic polynomials introduced in the previous section
reduces to estimates for character sums discussed in later sections.
For each fixed $\ell\geq1$, set
\[
z=N^{\delta/\ell}.
\]
Now, we define the truncation
\[
\omega_z(v)
=
\sum_{\substack{p\leq z\\p\nmid2aq}}
1_{p\mid v}.
\]
\begin{proposition}\label{prop:gaussian-moments}
Suppose that, uniformly in $n$ and $\eta$, we have
\[\frac1q\sum_{\substack{0 \leq r < q \\ r+nqK \notin S}}
\bigg\vert \sum_{\substack{p\leq z\\p\nmid2aq}}
\frac1p\left(\frac{\Delta_r}p\right) \bigg\vert^b
\ll_b (\log\log\log N)^b,
\qquad 1\leq b\leq 2\ell.\]
Then, uniformly in $n$ and $\eta$, and for all integer intervals $(I_r)_{r \in \{0,\dots,q-1\}} \subseteq \{0,1,\ldots,K\}$ with lengths at least $N^{1/4}$, we have
\begin{align*}
\frac{1}{\mathcal R_n} \sum_{\substack{0 \leq r < q \\ nqK + r \notin S}}
&\frac1{|I_r|}
\sum_{m\in I_r}
\bigl(\omega_z(f_r(m))
-\log\log N\bigr)^\ell \\
&=
\begin{cases}
(\ell-1)!!(\log\log N)^{\ell/2},&\ell\text{ even},\\
0,&\ell\text{ odd},
\end{cases}
+ O_{\ell}\left(
(\log\log\log N)
(\log\log N)^{(\ell-1)/2}
\right).
\end{align*}
\end{proposition}

It is a common strategy in proving results surrounding the Erd\H{o}s--Kac theorem to use generating functions (see, for example,
\cite[Section 2.3]{Lester2015}). We follow a similar approach in our proof.

\begin{proof}

We begin with the uncentered moments of $\omega_z(f_r(m))$ and consider its exponential generating function
\[
e^{t\omega_z(f_r(m))}
=
\sum_{j\geq0}
\frac{\omega_z(f_r(m))^j}{j!}t^j.
\]
Since $\omega_z(f_r(m)) = \sum_{p \leq z, \, p \nmid 2aq} 1_{p \mid f_r(m)}$ is a sum of indicator functions, and $e^{t1_A}=1+1_A(e^t-1)$, we have
\[
e^{t\omega_z(f_r(m))}
=
\prod_{\substack{p\leq z\\p\nmid2aq}}
\left(
1+1_{p\mid f_r(m)}
(e^t-1)
\right).
\]
Expanding the product and using that $1_{p \mid n} 1_{p' \mid n} = 1_{pp' \mid n}$ for distinct primes $p,p'$, and $n \in \N$, we obtain
\[
e^{t\omega_z(f_r(m))}
=
\sum_{\substack{
d\text{ square-free}\\
p\mid d\Rightarrow p\leq z,\ p\nmid2aq
}}
1_{d\mid f_r(m)}
(e^t-1)^{\omega(d)}.
\]
Averaged over $m \in I_r$, this is
\[\frac{1}{\vert I_r \vert} \sum_{m \in I_r} e^{t\omega_z(f_r(m))}
=
\sum_{\substack{
d\text{ square-free}\\
p\mid d\Rightarrow p\leq z,\ p\nmid2aq
}} (e^t-1)^{\omega(d)}
\frac{1}{\vert I_r \vert} \sum_{m \in I_r} 1_{d\mid f_r(m)}.\]
The second sum can be estimated by standard methods.
Since $p\nmid 2aq$, the quadratic polynomial $f_r$ has $1+\Legendre{\Delta_r}{p}$ roots modulo $p$. As $d$ is square-free, we apply the Chinese
remainder theorem and obtain
\[
\#\{m\bmod d:f_r(m)\equiv0\pmod d\}
=
\prod_{p\mid d}
\left(1+\Legendre{\Delta_r}{p}\right).
\]
Summing over these residue classes, we obtain
\[
\sum_{m\in I_r}1_{d\mid f_r(m)}
=
\frac{|I_r|}{d}
\prod_{p\mid d}
\left(1+\Legendre{\Delta_r}{p}\right)
+O\left(2^{\omega(d)}\right).
\]
For $p\nmid2aq$, define
\[
F_r(p)
:=
\frac1p\left(
1+\Legendre{\Delta_r}{p}
\right).
\]
 Then, dividing the previous identity by $|I_r|$, we obtain
\[\frac1{|I_r|}
\sum_{m\in I_r}
1_{d\mid f_r(m)}
=
\prod_{p\mid d}F_r(p)
+
O\left(\frac{2^{\omega(d)}}{|I_r|}\right).\]
Hence,
\[\frac1{|I_r|}
\sum_{m\in I_r}
e^{t\omega_z(f_r(m))}
=
\sum_{\substack{
d\text{ square-free}\\
p\mid d\Rightarrow p\leq z,\ p\nmid2aq
}}
(e^t-1)^{\omega(d)}
\left(\prod_{p\mid d}F_r(p)
+
O\left(\frac{2^{\omega(d)}}{|I_r|}\right)\right),
\]
where the main term is
\[\sum_{\substack{d \text{ square-free} \\ p \mid  d \Rightarrow p \leq z, \, p \nmid 2aq}} (e^t-1)^{\omega(d)} \prod_{p\mid d}F_r(p)
= \prod_{\substack{p\leq z\\p\nmid2aq}} \bigl(1+F_r(p)(e^t-1)\bigr).\]
For a power series $g$, let $[t^j]g(t)$ denote its
$j$-th coefficient. Then, for $1\leq j\leq2\ell$,
the centered moments are
\begin{align*}
\frac1{|I_r|}
\sum_{m\in I_r}
\bigl(\omega_z(f_r(m))-\log\log N\bigr)^j
&=
j![t^j]e^{-t\log\log N}
\frac1{|I_r|}
\sum_{m\in I_r}e^{t\omega_z(f_r(m))}\\
&=
\underbrace{
j![t^j]e^{-t\log\log N}
\prod_{\substack{p\leq z\\p\nmid2aq}}
\bigl(1+F_r(p)(e^t-1)\bigr)
}_{[1]}\\
&\quad+
\underbrace{
j![t^j]e^{-t\log\log N}
\sum_{\substack{
d\text{ square-free}\\
p\mid d\Rightarrow p\leq z,\ p\nmid2aq
}}
(e^t-1)^{\omega(d)}
O\left(\frac{2^{\omega(d)}}{|I_r|}\right)
}_{[2]}.
\end{align*}
Here the errors inside the sum are independent of $t$.
Let us first bound the error term $[2]$.
Since
\[
e^t-1=t+\sum_{k=2}^{\infty}\frac{t^k}{k!},
\]
the lowest degree of $(e^t-1)^k$ as a power series in $t$
is $k$. Therefore, only those $d$ that are products of
at most $i$ distinct primes contribute to the coefficient
of $t^i$.
For each number of prime factors $k\leq i$, there are
at most $z^k$ such $d$, and the coefficients of $t^i$
in $(e^t-1)^k$ are bounded in terms of $i$ and $k$.
Therefore,
\[
[t^i]
\sum_{\substack{
d\text{ square-free}\\
p\mid d\Rightarrow p\leq z,\ p\nmid2aq
}}
(e^t-1)^{\omega(d)}
O\left(\frac{2^{\omega(d)}}{|I_r|}\right)
=
O_i\left(
\frac1{|I_r|}\sum_{k=0}^{i}2^kz^k
\right)
=
O_i\left(\frac{z^i}{|I_r|}\right).
\]
Thus, using the Cauchy product formula, we obtain
\[
[2]
=
j!\sum_{i=0}^{j}
\frac{(-\log\log N)^{j-i}}{(j-i)!}
O_i\left(\frac{z^i}{|I_r|}\right)
=
O_\ell\left(
\frac{(1+\log\log N)^jz^j}{|I_r|}
\right).
\]
Since $z=N^{\delta/\ell}$, $|I_r|\geq N^{1/4}$,
$j\leq2\ell$, and $0<\delta<1/8$, this is bounded by
\[
O_\ell\left(
(\log\log N)^{2\ell}N^{2\delta-1/4}
\right)
=o(1).
\]
In particular, this estimate is uniform in all the parameters
and the choice of interval.
We now turn to the main term $[1]$. Let
\[
M_r(t)
=
e^{-t\log\log N}
\prod_{\substack{p\leq z\\p\nmid2aq}}
\bigl(1+F_r(p)(e^t-1)\bigr).
\]
We take the logarithm to turn the product defining $M_r(t)$
into a sum. Since $M_r(0)=1$, its formal logarithm is
well defined, and
\[
\log M_r(t)
=
-t\log\log N
+
\sum_{\substack{p\leq z\\p\nmid2aq}}
\log\bigl(1+F_r(p)(e^t-1)\bigr).
\]
Using
\[
\log(1+x)=\sum_{s\geq1}\frac{(-1)^{s+1}}s x^s,
\]
we obtain
\begin{align*}
\log M_r(t)
&=
-t\log\log N
+
\sum_{s\geq1}\frac{(-1)^{s+1}}s
(e^t-1)^s
\sum_{\substack{p\leq z\\p\nmid2aq}}F_r(p)^s\\
&=
-t\log\log N
+
(e^t-1)
\sum_{\substack{p\leq z\\p\nmid2aq}}
\left(
\frac1p+\frac1p \Legendre{\Delta_r}{p}
\right)+
\sum_{s\geq2}\frac{(-1)^{s+1}}s
(e^t-1)^s
\sum_{\substack{p\leq z\\p\nmid2aq}}F_r(p)^s.
\end{align*}
At each fixed coefficient order, only finitely many values
of $s$ contribute.
The middle term is
\[
(e^t-1)
\bigg(
\log\log N
+
O_{\ell}(\log\log\log N)
+
\sum_{\substack{p\leq z\\p\nmid2aq}}
\frac1p\left(\frac{\Delta_r}{p}\right)
\bigg).
\]
Indeed, by Mertens' theorem,
\[
\sum_{p\leq z}\frac1p
=
\log\log z+O(1)
=
\log\log N+\log(\delta/\ell)+O(1).
\]
Moreover, using the bound for reciprocal prime divisors
and $aq\ll N^{O(1)}$ from the major-arc setup, we have
\[
\sum_{\substack{p\leq z\\p\mid2aq}}\frac1p
\ll\log\log\log N.
\]
For every fixed $s\geq2$, the sum over primes in the last term
satisfies
\[
0 \leq \sum_{\substack{p\leq z\\p\nmid2aq}}F_r(p)^s
\leq 2^s\sum_p\frac1{p^s}
=O_s(1).
\]
In particular, the linear coefficient of $\log M_r(t)$ is
\[
[t]\log M_r(t)
=
\sum_{\substack{p\leq z\\p\nmid2aq}}F_r(p)-\log\log N
=
\sum_{\substack{p\leq z\\p\nmid2aq}}
\frac1p\left(\frac{\Delta_r}{p}\right)
+
O_{\ell}(\log\log\log N).
\]
Since $\mathcal R_n=q-\frac{|S_n|}{K}$ and $|S_n|/K=o(q)$, we have
$q/\mathcal R_n=1+o(1)$. Thus, the assumed estimate implies that, for $1\le b\le 2\ell$,
\[
\begin{aligned}
\frac1{\mathcal R_n}
\sum_{\substack{0\le r<q\\r\notin S}}
\left|[t]\log M_r(t)\right|^b
&\ll_{\ell}
\frac{q}{\mathcal R_n}
(\log\log\log N)^b\\
&\ll_{\ell}
(\log\log\log N)^b.
\end{aligned}
\]
Thus, write

\[\log M_r(t)=
(\log\log N)(e^t-1-t)
+
\bigl([t]\log M_r(t)\bigr)(e^t-1)+
\sum_{s\geq2}\frac{(-1)^{s+1}}s
(e^t-1)^s
\sum_{\substack{p\leq z\\p\nmid2aq}}F_r(p)^s.
\]
The last sum has uniformly bounded coefficients at every
fixed order and starts at degree two. We are left to recover the coefficients of $M_r(t)$.
Since
\[
M_r(t)
=
\exp(\log M_r(t))
=
\sum_{k\geq0}\frac{(\log M_r(t))^k}{k!},
\]
and $\log M_r(t)$ has zero constant term, only
$k\leq j$ contribute to the coefficient of $t^j$.
By the Cauchy product formula,
\[
j![t^j]M_r(t)
=
j!\sum_{k=1}^{j}\frac1{k!}
\sum_{\substack{i_1+\cdots+i_k=j\\
i_1,\ldots,i_k\geq1}}
\prod_{\nu=1}^{k}
[t^{i_\nu}]\log M_r(t).
\]
From the preceding decomposition, we have
\[
[t^i]\log M_r(t)
=
\begin{cases}
[t]\log M_r(t), & i=1,\\[2mm]
\displaystyle
\frac{\log\log N+[t]\log M_r(t)}{i!}
+O_i(1), & i\geq2,
\end{cases}
\]
where the error is uniform in all variables.

In the Cauchy product, each factor of
$\log\log N$ comes from a coefficient of
degree at least two, whereas a factor of
$[t]\log M_r(t)$ may come from a coefficient
of degree one. Consequently, a contribution
containing $b$ factors of $[t]\log M_r(t)$
can contain at most $\lfloor(j-b)/2\rfloor$
factors of $\log\log N$.

If $j$ is even, the unique contribution of maximal degree $j/2$ in $\log\log N$ that is independent of $[t]\log M_r(t)$ arises when $k=j/2$ and $i_1=\cdots=i_k=2$, taking the leading term from each coefficient.
This contributes
\[
\frac{j!}{(j/2)!}
\left(\frac{\log\log N}{2}\right)^{j/2}
=
(j-1)!!(\log\log N)^{j/2}.
\]
If $j$ is odd, no such contribution exists.
All remaining terms independent of
$[t]\log M_r(t)$ contain at most
$\lfloor(j-1)/2\rfloor$ factors of
$\log\log N$.

Only finitely many compositions
contribute, and the numerical coefficients arising
from their expansions are bounded in terms of $j$.
Therefore, we obtain by summing these contributions
for $1\leq j\leq2\ell$
\begin{align}
j![t^j]M_r(t)
={}&
\begin{cases}
(j-1)!!(\log\log N)^{j/2},&j\text{ even},\\
0,&j\text{ odd},
\end{cases}
\notag\\
&+
O_j\left(
(1+\log\log N)^{\lfloor(j-1)/2\rfloor}
+
\sum_{b=1}^{j}
(1+\log\log N)^{\lfloor(j-b)/2\rfloor}
\left|[t]\log M_r(t)\right|^b
\right).
\label{eq:model-moment}
\end{align}

We now average \eqref{eq:model-moment} over $r+nqK\notin S $.
The averaged error is bounded by
\[
O_{\ell}\left(
(1+\log\log N)^{\lfloor(j-1)/2\rfloor}
+
\sum_{b=1}^{j}
(1+\log\log N)^{\lfloor(j-b)/2\rfloor}
(\log\log\log N)^b
\right).
\]
Since
\[
\log\log\log N=o(\sqrt{\log\log N}),
\]
we have, for $1\leq b\leq j$,
\[
(1+\log\log N)^{\lfloor(j-b)/2\rfloor}
(\log\log\log N)^b
\ll_\ell
(\log\log\log N)(\log\log N)^{(j-1)/2}.
\]
Therefore, the averaged error is
\[
O_{\ell}\left(
(\log\log\log N)(\log\log N)^{(j-1)/2}
\right).
\]
Combining this with the estimate for $[2]$, we obtain,
for $1\leq j\leq2\ell$,
\begin{align*}
\frac1{\mathcal R_n}
\sum_{\substack{0 \leq r < q \\ r+nqK \notin S}}
\frac1{|I_r|}
\sum_{m\in I_r}
\bigl(\omega_z(f_r(m))-\log\log N\bigr)^j
={}&
\begin{cases}
(j-1)!!(\log\log N)^{j/2},&j\text{ even},\\
0,&j\text{ odd},
\end{cases}\\
&+
O_{\ell}\left(
(\log\log\log N)(\log\log N)^{(j-1)/2}
\right).
\end{align*}
We conclude the proof by taking $j=\ell$.
\end{proof}
In particular, taking in the above $j=2\ell$ we obtain
\[
\frac1{\mathcal R_n}
\sum_{\substack{0 \leq r < q \\ r+nqK \notin S}}
\frac1{|I_r|}
\sum_{m\in I_r}
|\omega_z(f_r(m))-\log\log N|^{2\ell}
\ll_{\ell}(\log\log N)^\ell.
\]
By H\"older's inequality, this also implies
\begin{equation}
\frac1{\mathcal R_n}
\sum_{\substack{0 \leq r < q \\ r+nqK \notin S}}
\frac1{|I_r|}
\sum_{m\in I_r}
|\omega_z(f_r(m))-\log\log N|^j
\ll_{\ell}(\log\log N)^{j/2},
\qquad 0\leq j\leq2\ell.
\label{eq:absolute-moments}
\end{equation}
These estimates will be used when removing the truncation
in Section 3.4.
\subsection{Average of Quadratic Character Sums}\label{subsec: legendre-symbol-sum}

In this section, we prove the following proposition, which provides the bounds needed in \Cref{prop:gaussian-moments}.

\begin{proposition}\label{prop:legendre-symbol-sum}
    For all ${\ell} \in \N$, we have that
    \[\frac{1}{q} \sum_{\substack{0 \leq r < q \\ nqK + r \notin S}} \bigg\vert \sum_{\substack{p \leq z \\ p \nmid 2aq}} \frac{1}{p} \Legendre{\Delta_r}{p} \bigg\vert^{\ell} \ll_{\ell} (\log \log \log N)^{\ell}.\]
\end{proposition}

Let us first introduce some additional notation. For fundamental discriminants $d$, we denote the Dirichlet character arising from the corresponding Kronecker symbol by $\chi_d(n) := \Legendre{d}{n}$. We also introduce the following notation.
\begin{definition}
    Recall that $\sqf(n)$ denotes the signed square-free part of $n \in \Z \setminus \{0\}$, and let
    \[
    g(n) :=
    \begin{cases}
        4\sqf(n), & \text{if } \sqf(n) \equiv 2,3 \pmod 4, \\
        \sqf(n), & \text{if } \sqf(n) \equiv 1 \pmod 4.
    \end{cases}
    \]
    In this section and only in this section, we also define $\sqf(0) := 1$ and for consistency $g(0):=1$.
\end{definition}

Note that $g(n)$ is a fundamental discriminant for all non-square $n \in \Z$ and therefore $\chi_{g(n)}$ is a primitive character. Moreover, $\chi_{g(n)}$ is non-principal whenever $g(n) \neq 1$, or equivalently, $n$ is not a square.

The following lemma allows us to work with these characters instead of the Legendre symbols.

\begin{lemma}\label{lem:switching Delta with g(Delta) in Legendre symbol sums}
    For any finite set of primes $\mathcal{P}$, any ${\ell} \in \N$, and any $\Delta \in \Z \setminus \{0\}$, we have
    \[\bigg\vert \sum_{p \in \mathcal{P}} \frac{1}{p} \Legendre{\Delta}{p} \bigg\vert^{\ell} \ll_{\ell} \bigg\vert \sum_{p \in \mathcal{P}} \frac{\chi_{g(\Delta)}(p)}{p} \bigg\vert^{\ell} + (\log \log \log \vert \Delta \vert)^{\ell}\]
\end{lemma}

\begin{proof}
    For all odd $p \nmid \Delta$ we have \[\Legendre{\Delta}{p} = \Legendre{g(\Delta)}{p} = \chi_{g(\Delta)}(p).\]
    Hence,
    \[\sum_{p \in \mathcal{P}} \frac{1}{p} \Legendre{\Delta}{p} = \sum_{p \in \mathcal{P}} \frac{\chi_{g(\Delta)}(p)}{p} + O \bigg( \frac{1}{2} + \sum_{p \mid \Delta, \, p \neq 2} \frac{1}{p} \bigg),\]
    where the sum in the error term is bounded by $\frac{31}{30}$ for $\vert \Delta \vert \leq e^e$ and by $\sum_{p \mid \vert \Delta \vert} \frac{1}{p} \ll \log \log \log \vert \Delta \vert$ for $\vert \Delta \vert > e^e$ due to Mertens' theorem. This proves the case ${\ell}=1$. The general statement follows as well using the standard estimate
    \[\vert x+y \vert^{\ell} \ll_{\ell} \vert x \vert^{\ell} + \vert y \vert^{\ell}\]
    for $x,y \in \R$.
\end{proof}

We now make use of the following newly discovered theorem.

\begin{theorem}[Quasi-Riemann Hypothesis {\cite[Theorem 1.1]{QuasiRH1112}}]\label{thm:quasiRH}
    Every Dirichlet $L$-function, including $\zeta(s)$, has no zeros for $\Re s > 11/12.$
\end{theorem}

\begin{remark}
At the time of writing, there are two proofs of the Quasi-Riemann Hypothesis: the one cited above for $\Re s > 11/12$ and \cite[Theorem 1.1]{QuasiRH78} for $\Re s > 7/8$, which also has a Lean certificate \cite{Lean-QuasiRH}. Since the latter implies the former and considering the classical review process has not been finalised, the authors chose to work with \cite{QuasiRH1112}. 
\end{remark}

The following lemma provides a convenient way to use \Cref{thm:quasiRH} in our setting.

\begin{lemma}[{\cite[Lemma 2.1]{Granville02}}]\label{lem:dirichletL}
    Let $s = \sigma + it$ with $|t| \le 3q$, and let $y \ge 2$ be a real number. Let $\frac{1}{2} \le \sigma_0 < \sigma$, and suppose that the rectangle $\{z : \sigma_0 < \operatorname{Re}(z) \le 1, \, |\operatorname{Im}(z) - t| \le y + 3\}$ contains no zeros of $L(z, \chi)$. Then
    $$|\log L(s, \chi)| \ll \frac{\log q}{\sigma - \sigma_0}.$$
    Further, putting $\sigma_1 = \min\left(\sigma_0 + \frac{1}{\log y}, \frac{\sigma + \sigma_0}{2}\right)$,
    $$\log L(s, \chi) = \sum_{n \le y} \frac{\Lambda(n)\chi(n)}{n^s \log n} + O\left(\frac{\log q}{(\sigma_1 - \sigma_0)^2} y^{\sigma_1 - \sigma}\right).$$
\end{lemma}

We remark that even though the formulation in \cite{Granville02} does not mention it explicitly, $\chi$ is a character modulo $q$ and should clearly be restricted to non-principal characters.

\begin{proof}[Proof of \Cref{prop:legendre-symbol-sum}]
    Let $y = \max \{(\log |g(\Delta_r)|)^{24},e^{24}+1 \}$ and note that, if $z \leq y$, Mertens' theorem implies 
    \[\left\vert \sum_{p \le z} \frac{1}{p}\Legendre{\Delta_r}{p} \right\vert \le \sum_{p \le z} \frac{1}{p}  \le \sum_{p \le y}\frac{1}{p} \ll \log\log\log \vert g(\Delta_r) \vert.\]
    Thus, we may assume $z > y$ from now on.
    By \Cref{lem:switching Delta with g(Delta) in Legendre symbol sums} and the convergence of $\sum_{p^k, \, k \geq 2} \frac{1}{k p^k}$,
    we have
    \[\left\vert \sum_{p \le z} \frac{1}{p}\Legendre{\Delta_r}{p} \right\vert \ll \left\vert \sum_{n \le z} \frac{\Lambda(n)\chi_{g(\Delta_r)}(n)}{n\log n} \right\vert + \log \log \log N,\]
    where $\chi_{g(\Delta_r)}$ is non-principal as $\Delta_r \neq \square$ by choice of $r$. Furthermore, $\Delta_r$ not being a square yields $y \geq 2$ as $\vert d \vert \geq 3$ for all fundamental discriminants $d$.
    Now let $\sigma_0 = 11/12$. Then \Cref{thm:quasiRH} implies that the rectangle $\{w : \sigma_0 < \Re w \le 1 : |\Im w| \le z+3\}$ contains no zeros of $L(w,\chi_{g(\Delta_r)})$. Moreover, we have $\sigma_0 + \frac{1}{\log z} < \sigma_0 + \frac{1}{\log y} < \frac{1+\sigma_0}{2}$. Thus, setting $s=1$, $\sigma_y = \sigma_0 + \frac{1}{\log y}$ and $\sigma_z = \sigma_0 + \frac{1}{\log z}$, we obtain from applying \Cref{lem:dirichletL} twice that
    \begin{align*}
        \sum_{n \le z}\frac{\Lambda(n)\chi_{g(\Delta_r)}(n)}{n\log n} &= \log L(1,\chi_{g(\Delta_r)}) + O\left(\frac{\log \vert g(\Delta_r) \vert}{(\sigma_z-\sigma_0)^2}z^{\sigma_z-1}\right) \\
        &= \sum_{n \le y}\frac{\Lambda(n)\chi_{g(\Delta_r)}(n)}{n\log n} + O\left(\frac{\log \vert g(\Delta_r) \vert}{(\sigma_y-\sigma_0)^2}y^{\sigma_y-1}\right) + O\left(\frac{\log \vert g(\Delta_r) \vert}{(\sigma_z-\sigma_0)^2}z^{\sigma_z-1}\right).
    \end{align*}
    Hence,
    \begin{equation}\label{eq:quasi-RH proof}
    \begin{aligned}
        \left\vert \sum_{p \le z} \frac{1}{p}\Legendre{\Delta_r}{p} \right\vert
        &\ll \left\vert \sum_{2 \le n \le y} \frac{\Lambda(n)\chi_{g(\Delta_r)}(n)}{n\log n} \right\vert + \frac{\log \vert g(\Delta_r) \vert}{(\sigma_y-\sigma_0)^2}y^{\sigma_y-1} + \frac{\log \vert g(\Delta_r) \vert}{(\sigma_z-\sigma_0)^2}z^{\sigma_z-1} + \log\log\log N \\
        &\ll \sum_{p\le y}\sum_{k=1}^\infty \frac{1}{kp^k} + \frac{(\log y)^2}{\log \vert g(\Delta_r) \vert} + \frac{(\log z)^2\log \vert g(\Delta_r) \vert}{z^{1/12}} + \log \log \log N \\
        &\ll \log \log \log N
    \end{aligned}
    \end{equation}
    using $\vert g(\Delta_r) \vert \ll \vert \Delta_r \vert \ll N^2$. Finally, as
    \[\bigg\vert \sum_{p \mid 2aq} \frac{1}{p} \Legendre{\Delta_r}{p} \bigg\vert \ll \log \log \log q \ll \log \log \log N,\]
    \Cref{eq:quasi-RH proof} proves \Cref{prop:legendre-symbol-sum}.
\end{proof}
\subsection{Exceptional Primes and Truncation}

In the prior sections, we established the required moment estimates for the more convenient $\omega_z(f_r(m)).$ In the following, we return to the full prime-counting function $\omega(f_r(m))$ and show that the contribution of the exceptional primes $p \mid aq$ is sufficiently small. Finally, we combine  all previous results, bound the contribution of primes removed by the truncation, and provide the proof of \Cref{prop:Major Arcs}, resolving the major-arc case of \Cref{thm:moment theorem}.

\begin{lemma}\label{lem:generalized gcd condition lemma}
    Let $(I_r)_{r \in \{0,\dots,q-1\}}$ be intervals of integers of lengths $\vert I_r \vert \geq N^{1/4}$, and let $b \le 2\ell$. Then
        $$\frac{1}{q}\sum_{0\le r < q} \frac{1}{|I_r|} \sum_{m\in I_r} \bigg(\sum_{p\mid aq} 1_{p \mid f_{n,r,\eta_r}(m)}\bigg)^b = O_\ell((\log\log\log N)^b),$$
    where whenever we sum over primes, we assume $2 < p \leq z$. 
\end{lemma}
\begin{proof}
    We first define $\stirlingii{b}{k}$ as the Stirling numbers of the second kind, and set $C_{k,j} = \stirlingii{b}{k}\binom{k}{j}j!(k-j)! = \stirlingii{b}{k}k!$. We rewrite the sum as 
    \begin{align*}
        &\quad \frac{1}{q}\sum_{0\le r < q} \frac{1}{|I_r|} \sum_{m\in I_r} \sum_{k=1}^b\stirlingii{b}{k}\sum_{\substack{p_1,\ldots,p_k\mid aq \\ p_i \text{ distinct}}} 1_{\prod_{i=1}^{k}p_i \mid f_{n,r,\eta_r}(m)} \\
        &= \frac{1}{q}\sum_{k=1}^b\stirlingii{b}{k}\sum_{0\le r < q} \frac{1}{|I_r|} \sum_{m\in I_r} \sum_{\substack{p_1,\ldots,p_k\mid aq \\ p_i \text{ distinct}}} 1_{\prod_{i=1}^{k}p_i \mid f_{n,r,\eta_r}(m)} \\
        &\le  \frac{1}{q}  \sum_{k=1}^b\stirlingii{b}{k} \sum_{\eta \in \{0,\pm1  \}}\sum_{0\le r < q}\frac{1}{|I_r|} \sum_{m\in I_r}\sum_{\substack{p_1,\ldots,p_k\mid aq \\ p_i \text{ distinct}}}1_{\prod_{i=1}^{k}p_i \mid f_{n,r,\eta}(m)} \\
        &\le \max_{\eta\in \{0,\pm1\}}3 \frac{1}{q}  \sum_{k=1}^b\stirlingii{b}{k}\sum_{0\le r < q}\frac{1}{|I_r|} \sum_{m\in I_r} \sum_{\substack{p_1,\ldots,p_k\mid aq \\ p_i \text{ distinct}}}1_{\prod_{i=1}^{k}p_i \mid f_{n,r,\eta}(m)}.
    \end{align*}
    But since $(a,q)=1$, we have
    \begin{align*}
        &\le 3\max_{\eta\in \{0,\pm1\}} \frac{1}{q}  \sum_{k=1}^b\sum_{j=0}^kC_{k,j}\sum_{0\le r < q}\frac{1}{|I_r|} \sum_{m\in I_r} \sum_{\substack{p_1,\ldots,p_j\mid q \\ p_{j+1},\ldots,p_{k} \mid a \\ p_i \text{ distinct}}}1_{\prod_{i=1}^{k}p_i \mid f_{n,r,\eta}(m)} \\
        &\le 3\max_{\eta\in \{0,\pm1\}}\sum_{k=1}^b\sum_{j=0}^kC_{k,j}\frac{1}{q}\sum_{0\le r < q} \frac{1}{|I_r|} \sum_{m\in I_r}\sum_{\substack{P_j'\mid q \\ P_{k-j} \mid a}}1_{P_j'P_{k-j} \mid f_{n,r,\eta}(m)}, \\
        &\le 3\max_{\eta\in \{0,\pm1\}}\sum_{k=1}^b\sum_{j=0}^kC_{k,j} \frac{1}{q}  \sum_{0\le r < q} \frac{1}{|I_r|} \sum_{m\in I_r}\left(\sum_{P_j\mid q }1_{P_j \mid f_{n,r,\eta}(m)}\right)\left(\sum_{P_{k-j}\mid a}1_{P_{k-j} \mid f_{n,r,\eta}(m)}\right)
    \end{align*}
    where $P_i$ ranges over ordered products of $i$ distinct primes $p\in (2,z]$. We then see that the sum over $P_{k-j}\mid a$ is independent of $m$, and thus we can write our sum as 
    $$3\max_{\eta\in \{0,\pm1\}}\sum_{k=1}^b\sum_{j=0}^kC_{k,j} \frac{1}{q}  \sum_{0\le r < q} \sum_{P_{k-j}\mid a}1_{P_{k-j} \mid f_{n,r,\eta}(0)}\frac{1}{|I_r|} \sum_{m\in I_r} \sum_{P_j\mid q }1_{P_j \mid f_{n,r,\eta}(m)},$$
    for which the inner sum is
    \begin{align*}
        \sum_{P_j \mid q}\frac{1}{|I_r|} \sum_{m\in I_r}1_{P_j \mid 2arm + \lfloor \alpha(nqK)^2 + \frac{ar^2}{q}\rfloor + 2nKar + \eta} &\le
    \sum_{P_j \mid q}\left(\frac{(r,P_j)}{P_j} + O\left(\frac{P_j}{|I_r|}\right)\right) \\
    &= \sum_{P_j \mid q}\left(\frac{1}{P_j}\sum_{\substack{d \mid P_j \\ d \mid r}}\varphi(d) + O\left(\frac{P_j}{|I_r|}\right)\right).
    \end{align*}
    We can bring the error term outside all the sums as $|I_r| \ge N^{1/4} \gg_b \omega(aq)^b.$ Thus omitting the $O(1)$ error term, we bound the entire sum by  
    $$3\max_{\eta\in \{0,\pm1\}}\sum_{k=1}^b\sum_{j=0}^kC_{k,j}\sum_{P_j \mid q}\frac{1}{P_j}\sum_{d \mid P_j}\varphi(d)\sum_{P_{k-j}\mid a}\frac{1}{q}\sum_{0 \le r < q}1_{d\mid r}1_{P_{k-j}\mid f_{n,r,\eta}(0)},$$
    which, when writing $ds = r$ becomes 
    \begin{equation}\label{eq: exceptional primes pre Erdos-Turan}
        3\max_{\eta\in \{0,\pm1\}}\sum_{k=1}^b\sum_{j=0}^kC_{k,j}\sum_{P_j \mid q}\frac{1}{P_j}\sum_{d \mid P_j}\varphi(d)\sum_{P_{k-j}\mid a}\frac{1}{q}\sum_{0 \le s <\frac{q}{d}}1_{P_{k-j}\mid f_{n,ds,\eta}(0)}.
    \end{equation}
    Here, we use that $k\mid\lfloor \gamma \rfloor \iff \{\gamma / k\} < 1/k,$ combined with \Cref{lem:Erdos-Turan} inequality with $M=\lfloor\sqrt{q}\rfloor$ to obtain
    \begin{align*}
    &\quad \sum_{P_{k-j}\mid a}\frac{1}{q}\sum_{0 \le s <\frac{q}{d}}1_{P_{k-j}\mid f_{n,ds,\eta}(0)} \\
    &\le \frac{1}{q}  \sum_{P_{k-j} \mid a} \#\left\{0 \le s < \frac{q}{d} : \left\{\frac{\alpha q(nqK)^2 + a(ds)^2 + 2nKaqds - \eta q}{qP_{k-j}}\right\} < \frac{1}{P_{k-j}}\right\} \\     
    &\le\frac{1}{q}  \sum_{P_{k-j}\mid a}\left( \frac{q/d}{P_{k-j}} + O\left(\frac{q}{M} + \sum_{k=1}^{M}\frac{1}{k}\left|\sum_{0\le s < \frac{q}{d}}e\left(k\frac{\alpha q(nqK)^2 + a(ds)^2 + 2nKaqds - \eta q}{qP_{k-j}}\right)\right|\right) \right)\\
    &\le\frac{1}{q}  \sum_{P_{k-j}\mid a}\left( \frac{q/d}{P_{k-j}} + O\left(\sqrt{q} + \sum_{k\le\sqrt{q}}\frac{1}{k}\left|\sum_{0\le s < \frac{q}{d}}e\left(\frac{2nKaqkds + ak(ds)^2}{qP_{k-j}}\right)\right|\right) \right) \\
    &=\frac{1}{q}  \sum_{P_{k-j}\mid a}\left( \frac{q/d}{P_{k-j}} + O\left(\sqrt{q} + \sum_{k\le\sqrt{q}}\frac{1}{k}\left|\sum_{0\le s < \frac{q}{d}}e\left(\frac{(a/P_{k-j})k(ds)^2}{q}\right)\right|\right)\right) \\
    &=\frac{1}{q}  \sum_{P_{k-j}\mid a}\left( \frac{q/d}{P_{k-j}} + O\left(\sqrt{q} + \sum_{k\le\sqrt{q}}\frac{1}{k}\left|\sum_{0\le s < \frac{q}{d}}e\left(\frac{(a/P_{k-j})kds^2}{q/d}\right)\right|\right)\right)
    \end{align*}
    The inner quadratic Gauss sum is bounded by
    \begin{align*}
        &\le \frac{1}{q}  \sum_{P_{k-j}\mid a} \left(\frac{q/d}{P_{k-j}} + O\left(\sqrt{q} + \sum_{k\le\sqrt{q}}\frac{\sqrt{\frac{q}{d}\cdot (kd,q/d)}}{k}\right)\right) \\
        &\le \frac{1}{q}  \sum_{P_{k-j}\mid a} \left(\frac{q/d}{P_{k-j}} + O\left(\sqrt{q} + \sum_{k\le\sqrt{q}}\frac{\sqrt{q\cdot (k,q/d)}}{k}\right)\right) \\
        &\le\frac{1}{d} \sum_{P_{k-j}\mid a} \frac{1}{P_{k-j}} + O_\varepsilon\left(\binom{\omega(a)}{k-j}\left(q^{-1/2} + q^{-1/2 + \varepsilon}\log q \right)\right) \\
        &= \frac{1}{d}\sum_{P_{k-j}\mid a} \frac{1}{P_{k-j}} + O_\varepsilon\left(q^{-1/2+2\varepsilon}\right),
    \end{align*}   
    where $\binom{\omega(a)}{k-j}\log q \ll \log(a)^{k-j}\log q \ll q^\varepsilon$. We can then plug this back into \Cref{eq: exceptional primes pre Erdos-Turan} and get a main term of 
    $$3\sum_{k=1}^b\sum_{j=0}^kC_{k,j}\sum_{P_{k-j}\mid a}\frac{1}{P_{k-j}}\sum_{P_j \mid q}\frac{1}{P_j}\sum_{d \mid P_j}\frac{\varphi(d)}{d}$$
    and an error term of 
    $$O_{b,\varepsilon}\left(3\sum_{k=1}^b\sum_{j=0}^kC_{k,j}\sum_{P_j \mid q}\frac{1}{P_j}\sum_{d \mid P_j}\varphi(d)q^{-1/2 + 2\varepsilon}\right).$$
    Since the inner sums of our main term are multiplicative, and $P_j$ is square-free, we have that 
    $$\sum_{P_j\mid q} \frac{1}{P_j}\sum_{d \mid P_j} \frac{\varphi(d)}{d} = \sum_{P_j\mid q} \frac{1}{P_j}\prod_{p \mid P_j} \left(\frac{\varphi(p)}{p} + \frac{\varphi(1)}{1}\right) =\sum_{P_j\mid q} \frac{1}{P_j}\prod_{p \mid P_j} \left(2 - \frac{1}{p}\right) \le 2^j\sum_{P_j\mid q} \frac{1}{P_j} \ll_j (\log\log\log q)^j.$$
    Furthermore, we know that $\sum_{P_{k-j}\mid a}\frac{1}{P_{k-j}} \ll (\log\log\log a)^{k-j}$ and using $\sum_{d\mid P_j}\varphi(d) = P_j$, we have that the sum is bounded by 
    \begin{align*}
        &\ll \sum_{k=1}^b\sum_{j=0}^kC_{k,j}\left((\log\log\log q)^j(\log\log\log a)^{k-j} + q^{-1/2+\varepsilon}\sum_{P_j \mid q}\frac{1}{P_j}\sum_{d \mid P_j}\varphi(d) \right)\\
        &\ll \sum_{k=1}^b\sum_{j=0}^kC_{k,j} \left((\log\log\log N)^k + q^{-1/2+\varepsilon}\binom{\omega(q)}{j}\right)\\ 
        &\ll \sum_{k=1}^b\sum_{j=0}^kC_{k,j} (\log\log\log N)^k \\ 
        &\ll_b (\log\log\log N)^b,
    \end{align*}
    concluding our proof that 
    \[\frac{1}{q}\sum_{0\le r < q} \frac{1}{|I_r|} \sum_{m\in I_r} \left(\sum_{p\mid aq} 1_{p \mid f_{n,r,\eta_r}(m)}\right)^b = O_\ell((\log\log\log N)^b).\]
\end{proof}

We are now able to provide the proof of \Cref{prop:Major Arcs}.

\begin{proof}[Proof of \Cref{prop:Major Arcs}]
    We first show that
    \begin{equation}\label{eq:truncation final proof 1}
    \begin{multlined}
    \frac{1}{\mathcal{R}_n} \sum_{\substack{0 \leq r < q \\ nqK+r \notin S}} \frac{1}{\vert I_r \vert} \sum_{m \in I_r} (\omega(f_{n,r,\eta_r}(m)) - \log \log N)^{\ell}
    = \begin{cases}
        ({\ell}-1)!! (\log \log N)^{{\ell}/2}, & {\ell} \text{ even}, \\
        0, & {\ell} \text{ odd.}
    \end{cases} \\
    \qquad \qquad \qquad + O_{\ell}((\log \log \log N) (\log \log N)^{({\ell}-1)/2})
    \end{multlined}
    \end{equation}
    for all integer intervals $(I_r)_{r \in \{0, \dots, q-1\}}$ with lengths at least $N^{1/4}$ and uniformly in $n$ and $\eta$.
    In order to deduce this from \Cref{prop:gaussian-moments}, which uses \Cref{prop:legendre-symbol-sum}, we need to transition from $\omega$ to $\omega_z$. Therefore, let $\mathcal{E}$ denote the set of primes $p$ with $p \mid 2aq$ or $p > z$ and observe that
    \begin{align*}
        \frac{1}{\mathcal{R}_n} \sum_{\substack{0 \leq r < q \\ nqK+r \notin S}} &\frac{1}{\vert I_r \vert} \sum_{m \in I_r} (\omega(f_{n,r,\eta_r}(m)) - \log \log N)^{\ell} \\
        &= \frac{1}{\mathcal{R}_n} \sum_{\substack{0 \leq r < q \\ nqK+r \notin S}} \frac{1}{\vert I_r \vert} \sum_{m \in I_r} \bigg( \omega_z(f_{n,r,\eta_r}(m)) + \sum_{p \in \mathcal{E}} 1_{p \mid f_{n,r,\eta_r}(m)} - \log \log N \bigg)^{\ell} \\
        &= \frac{1}{\mathcal{R}_n} \sum_{\substack{0 \leq r < q \\ nqK+r \notin S}} \frac{1}{\vert I_r \vert} \sum_{m \in I_r} \sum_{\nu = 0}^{\ell} \binom{\ell}{\nu} (\omega_z(f_{n,r,\eta_r}(m)) - \log \log N)^\nu \bigg( \sum_{p \in \mathcal{E}} 1_{p \mid f_{n,r,\eta_r}(m)} \bigg)^{{\ell}-\nu} \\
        &= \sum_{\nu = 0}^{\ell} \binom{\ell}{\nu} \frac{1}{\mathcal{R}_n} \sum_{\substack{0 \leq r < q \\ nqK+r \notin S}} \frac{1}{\vert I_r \vert} \sum_{m \in I_r} (\omega_z(f_{n,r,\eta_r}(m)) - \log \log N)^\nu \bigg(\sum_{p \in \mathcal{E}} 1_{p \mid f_{n,r,\eta_r}(m)} \bigg)^{{\ell}-\nu}.
    \end{align*}
    For $\nu = {\ell}$, \Cref{prop:gaussian-moments} and \Cref{prop:legendre-symbol-sum} now yield the correct main term
    \[\begin{cases}
        ({\ell}-1)!! (\log \log N)^{{\ell}/2}, & {\ell} \text{ even}, \\
        0, & {\ell} \text{ odd.}
    \end{cases}\]
    with an error of $O_{\ell}((\log \log \log N)(\log \log N)^{({\ell}-1)/2})$. 
    For $0 \leq \nu < {\ell}$, we apply the Cauchy-Schwarz inequality and obtain
    \begin{align*}
        \bigg\vert \frac{1}{\mathcal{R}_n} \sum_{\substack{0 \leq r < q \\ nqK+r \notin S}} \frac{1}{\vert I_r \vert} \sum_{m \in I_r} (\omega_z(f_{n,r,\eta_r}(m)) &- \log \log N)^\nu \bigg(\sum_{p \in \mathcal{E}} 1_{p \mid f_{n,r,\eta_r}(m)} \bigg)^{{\ell}-\nu} \bigg\vert \\
        &\ \leq \underbrace{\bigg( \frac{1}{\mathcal{R}_n} \sum_{\substack{0 \leq r < q \\ nqK+r \notin S}} \frac{1}{\vert I_r \vert} \sum_{m \in I_r} (\omega_z(f_{n,r,\eta_r}(m)) - \log \log N)^{2\nu} \bigg)^{1/2}}_{[1]} \\
        &\quad \quad \cdot \underbrace{\bigg( \frac{1}{\mathcal{R}_n} \sum_{\substack{0 \leq r < q \\ nqK+r \notin S}} \frac{1}{\vert I_r \vert} \sum_{m \in I_r} \bigg(\sum_{p \in \mathcal{E}} 1_{p \mid f_{n,r,\eta_r}(m)} \bigg)^{2({\ell}-\nu)} \bigg)^{1/2}}_{[2]}.
    \end{align*}
    The first factor [1] can be bounded by another application of \Cref{prop:gaussian-moments} in the form of \Cref{eq:absolute-moments} and \Cref{prop:legendre-symbol-sum}, resulting in
    \begin{align*}
        \bigg( \frac{1}{\mathcal{R}_n} \sum_{\substack{0 \leq r < q \\ nqK+r \notin S}} \frac{1}{\vert I_r \vert} \sum_{m \in I_r} (\omega_z(f_{n,r,\eta_r}(m)) - \log \log N)^{2\nu} \bigg)^{1/2}
        &\ll_{\nu} (\log \log N)^{\nu/2}.
    \end{align*}
    For the second factor [2], we cover $\mathcal{E}$ by the set of primes in $(2,z]$ dividing $aq$, the set of primes that are larger than $z$, and $\{2\}$. The standard estimate $\vert x+y \vert^m \ll_m \vert x \vert^m + \vert y \vert^m$ for $x,y \in \R$ then implies
    \begin{align*}
        \frac{1}{\mathcal{R}_n} \sum_{\substack{0 \leq r < q \\ nqK+r \notin S}} \frac{1}{\vert I_r \vert} \sum_{m \in I_r} \bigg(\sum_{p \in \mathcal{E}} 1_{p \mid f_{n,r,\eta_r}(m)} \bigg)^{2({\ell}-\nu)}
        &\ll_{\ell} \frac{1}{\mathcal{R}_n} \sum_{\substack{0 \leq r < q \\ nqK+r \notin S}} \frac{1}{\vert I_r \vert} \sum_{m \in I_r} \bigg(\sum_{\substack{p \mid aq \\ 2 < p \leq z}} 1_{p \mid f_{n,r,\eta_r}(m)} \bigg)^{2({\ell}-\nu)} \\
        &+ \frac{1}{\mathcal{R}_n} \sum_{\substack{0 \leq r < q \\ nqK+r \notin S}} \frac{1}{\vert I_r \vert} \sum_{m \in I_r} \bigg(\sum_{p > z} 1_{p \mid f_{n,r,\eta_r}(m)} \bigg)^{2({\ell}-\nu)} + O(1).
    \end{align*}
    The left summand is $O_{\ell}((\log \log \log N)^{2({\ell}-\nu)})$ by \Cref{lem:generalized gcd condition lemma}, and since there are $\ll \log_z N$ primes $p > z$ dividing $f_{n,r,\eta_r}(m)$, the right summand satisfies
    \[\frac{1}{\mathcal{R}_n} \sum_{\substack{0 \leq r < q \\ nqK+r \notin S}} \frac{1}{\vert I_r \vert} \sum_{m \in I_r} \bigg(\sum_{p > z} 1_{p \mid f_{n,r,\eta_r}(m)} \bigg)^{2({\ell}-\nu)} 
    \ll_{\ell} \frac{1}{\mathcal{R}_n} \sum_{\substack{0 \leq r < q \\ nqK+r \notin S}} \frac{1}{\vert I_r \vert} \sum_{m \in I_r} \bigg( \frac{\log N}{\log z} \bigg)^{2({\ell}-\nu)} \ll_{\ell} 1.\]
    Combining all these estimates, we see that each $\nu < {\ell}$ contributes
    \[O_{\ell}((\log \log \log N)^{{\ell}-\nu} (\log \log N)^{\nu/2}),\]
    proving \Cref{eq:truncation final proof 1}.
    We are now ready to substitute this in the bounds from \Cref{subsec: alpha to polynomials} and \Cref{eq: double bounds}. This gives
\begin{multline*}
    \min_{\substack{\eta\in\{0,\pm1\}^q \\ 0 \le n \le \frac{N}{qK}-1}}\begin{cases}
        ({\ell}-1)!! (\log \log N)^{{\ell}/2}, & {\ell} \text{ even}, \\
        0, & {\ell} \text{ odd.}
    \end{cases}+O_{\ell}((\log \log \log N) (\log \log N)^{({\ell}-1)/2}) \\
    \le\frac{1}{|\mathcal C_N|}\sum_{n\in \mathcal C_N} \left(\omega(\lfloor\alpha n^2\rfloor)-\log\log N \right)^\ell \leq \\
    \max_{\substack{\eta\in\{0,\pm1\}^q \\ 0 \le n \le \frac{N}{qK}-1}}\begin{cases}
        ({\ell}-1)!! (\log \log N)^{{\ell}/2}, & {\ell} \text{ even}, \\
        0, & {\ell} \text{ odd.}
    \end{cases}+O_{\ell}((\log \log \log N) (\log \log N)^{({\ell}-1)/2}).
\end{multline*}
    Since the sum is bounded from both sides, we have
    \begin{align*}
        \frac{1}{|\mathcal C_N|}\sum_{n\in \mathcal C_N} \left(\omega(\lfloor\alpha n^2\rfloor)-\log\log N \right)^\ell &= \begin{cases}
        ({\ell}-1)!! (\log \log N)^{{\ell}/2}, & {\ell} \text{ even}, \\
        0, & {\ell} \text{ odd},
        \end{cases}\\
        &\qquad \qquad \qquad+O_{\ell}((\log \log \log N) (\log \log N)^{({\ell}-1)/2}),
    \end{align*}
    completing the proof of \Cref{prop:Major Arcs}.
\end{proof}

\section{The Rational Case}
In this section, we will focus on the rational case. This is much simpler to deal with, as for $\alpha = \frac{a}{q}$, we split $\lfloor\alpha n^2\rfloor$ into $q$ different integer polynomials and are only concerned with their irreducibility. 

We start by proving \Cref{thm:rational theorem}. That is,  \textit{the number of distinct prime factors of $\lfloor \alpha n^2\rfloor$ obeys a Gaussian law if and only if $\sqf(a) > q$.} The proof relies on a result by Minoru Tanaka about integer polynomials. 

\begin{theorem}[\cite{Tanaka55}]\label{thm:tanaka}
    Let $f \in \Z[X]$ be a polynomial with $k$ distinct irreducible factors. The number of distinct prime factors of $f(n)$ follows a Gaussian distribution of mean and variance $k\log\log n$. 
\end{theorem}
\begin{proof}[Proof of \Cref{thm:rational theorem}]
    We write $a = a_1a_2^2$, where $a_1 = \sqf(a)$, and $a_2 = \sqp(a)$. We convert $\lfloor \alpha n^2\rfloor$ to a polynomial. For a fixed $n \in \N$, we can write $n = qm + r$, $0\le r<q$.  This gives
    $$\lfloor \alpha n^2\rfloor = \left\lfloor \frac{a}{q}(qm+r)^2\right\rfloor = aqm^2 + 2arm + \lfloor\frac{ar^2}{q}\rfloor.$$
    We thus get $q$ polynomials. For $r \in \{0,1,\ldots,q-1\}$, we denote $f_r(m) = aqm^2 + 2arm + \lfloor\frac{ar^2}{q}\rfloor$. The discriminant $\Delta_r$ of $f_r$ is 
    $$\Delta_r = 4aq\left\{\frac{ar^2}{q}\right\} = 4a(ar^2 \bmod q) = (2a_2)^2a_1(ar^2\bmod q).$$
    As $a_1 > q$ and is square-free, we know that any value $x$ $\bmod q$ will not allow $a_1x$ to be a non-zero square, and thus $f_r$ has at most one irreducible factor. We conclude by applying \Cref{thm:tanaka}. 

    For the forward direction, we know that $a_1 < q$. We will show that the average is not $\log\log n$. As for the backward direction, we write 
    $f_r(m) = aqm^2 + 2arm + \lfloor\frac{ar^2}{q}\rfloor$. The discriminant thus can be written as
    $$\Delta_r =  4a(ar^2 \bmod q) = 4a\left(a_1(a_2r)^2\bmod q\right).$$
    Since $(a_2,q)=1$, there is an $r$ such that $a_2r \equiv 1 \bmod q$. For such an $r$, we have $\Delta_r = 4a(a_1 \bmod q) = (2a_1a_2)^2 > 0$, meaning $f_r$ has two distinct factors. Let us define a set $S$ of $r$'s yielding $\Delta_r$ being a strictly positive square.
    From \Cref{thm:tanaka}, we have that $\omega(\lfloor \alpha n^2 \rfloor)$ behaves like the sum of two Gaussian, one with mean and variance $\log\log n$ that has weight $\frac{q-|S|}{q}$ and another one with mean and variance $2\log\log n$ that has weight $\frac{|S|}{q}$.  
\end{proof}

To complete the argument, we give a brief overview of the size of $|S|$ through the proof of \Cref{cor: rational corollary}. 
\begin{proof}
    For $f_r$ to have two distinct irreducible factors, we need $\Delta_r = (2a_2)^2a_1\left(a_1(a_2r)^2\bmod q\right)$ to be a perfect square. Remember that $a_2(\Z/q\Z)\cong \Z/q\Z$. In particular, the number of values of $r$ for which $\Delta_r$ is a positive perfect square is equivalent to the number of $r$ with $a_1(a_1r^2 \bmod q)$ being a perfect square. As $(a_1,q) = 1$ and $a_1$ is square-free, we need $(a_1r^2 \bmod q) = k^2a_1$, $k\in\N$. As we are working modulo $ q$, only the smallest $$\left\lfloor\sqrt{\frac{q-1}{a_1}}\right\rfloor$$ integers $k$ are suitable. 

    This means that the set $S$ from the proof of \Cref{thm:rational theorem} has size 
    $$\sum_{k\le \left\lfloor\sqrt{\frac{q-1}{a_1}}\right\rfloor}\#\{0\le r<q:r^2\equiv k^2\bmod q\}.$$
    Now, for every $r \in S$, we know that the polynomial $f_r(m) = \lfloor \alpha (mq+r)^2\rfloor$ will have two distinct irreducible factors, and only one when $r\notin S$. Thus, for any $\tau\in\R$, by \Cref{thm:tanaka}, we have 
    $$\lim_{N\to\infty} \frac{1}{N} \#\left\{n \le N : \frac{\omega(\lfloor \alpha n^2 \rfloor)-(1 + 1_{S}(n))\log\log N}{\sqrt{(1+1_{S}(n))\log\log N}} < \tau\right\} = \Phi(\tau).$$
\end{proof}
\section{Final Remarks}
We close by discussing possible generalizations of our results. Having shown that, for positive irrational $\alpha$, $\omega(\lfloor \alpha n^2\rfloor)$ obeys a Gaussian distribution, one may ask whether the same holds for arbitrary quadratic polynomials with at least one non-constant irrational coefficient. This leads us to the following conjecture.

\begin{conjecture}
Let $\alpha > 0,\beta,\gamma$ be real numbers, with at least one of $\alpha$ and $\beta$ irrational. The number of distinct prime factors of $\lfloor \alpha n^2+\beta n+\gamma\rfloor$ obeys a Gaussian law. More precisely, we have
    $$\lim_{N\to\infty} \frac{1}{N} \#\left\{n \le N : \frac{\omega(\lfloor \alpha n^2 +\beta n+\gamma \rfloor)-\log\log N}{\sqrt{\log\log N}} < \tau\right\} = \Phi(\tau).$$
\end{conjecture}
We also expect a similar result to hold for higher-degree polynomials.
\begin{conjecture}\label{conj: higher degrees}
Let $\alpha_k > 0,\alpha_{k-1},\ldots,\alpha_0$ be real numbers, with at least one of $\alpha_k,\alpha_{k-1},\ldots,\alpha_1$ irrational. The number of distinct prime factors of $\lfloor \alpha_k n^k+\alpha_{k-1}n^{k-1}+\cdots+\alpha_0\rfloor$ obeys a Gaussian law. More precisely, we have
    $$\lim_{N\to\infty} \frac{1}{N} \#\left\{n \le N : \frac{\omega(\lfloor \alpha_k n^k+\alpha_{k-1}n^{k-1}+\cdots+\alpha_0 \rfloor)-\log\log N}{\sqrt{\log\log N}} < \tau\right\} = \Phi(\tau).$$
\end{conjecture}
In fact, adapting the techniques developed in \Cref{sec: Minor arcs}, one may prove that \Cref{conj: higher degrees} holds when the leading coefficient is irrational and of finite Diophantine type.

\newpage
\printbibliography

\bigskip
\footnotesize
\noindent
Catinca (Cati) Alexandru\\
\textsc{Brandeis University}\\
\href{mailto:catialexandru@brandeis.edu}
{\texttt{catialexandru@brandeis.edu}}

\bigskip
\footnotesize
\noindent
Johannes Borgwardt\\
\textsc{Georg-August-Universit\"{a}t G\"{o}ttingen}\\
\href{mailto:johannes.borgwardt@stud.uni-goettingen.de}
{\texttt{johannes.borgwardt@stud.uni-goettingen.de}}

\bigskip
\footnotesize
\noindent
Benjamin Faltin\\
\textsc{{\'E}cole Polytechnique F{\'e}d{\'e}rale de Lausanne (EPFL)}\\
\href{mailto:benjamin.faltin@epfl.ch}
{\texttt{benjamin.faltin@epfl.ch}}

\bigskip
\footnotesize
\noindent
Giorgos Petritakis\\
\textsc{University of Crete}\\
\href{mailto:giorgospetritakis@gmail.com}
{\texttt{giorgospetritakis@gmail.com}}
\end{document}